\documentclass[11pt]{amsart}
\usepackage[a4paper,margin=28mm]{geometry}

\usepackage{xcolor}
\usepackage{lineno}
\usepackage{amsmath,amssymb,mathrsfs}
\usepackage[shortlabels]{enumitem}
\usepackage{microtype}

\usepackage[
    pagebackref=true,
    colorlinks=false,
    pdfborder={0 0 1},
    linkbordercolor={1 0 0},
    citebordercolor={0 1 0},
    urlbordercolor={0 1 1}
]{hyperref}

\newcommand{\dbar}{\bar\partial}
\newcommand{\step}[1]{\par\smallskip\noindent\textbf{Step #1.}\ }
\newcommand{\eps}{\varepsilon}
\newcommand{\C}{\mathbb C}

\DeclareMathOperator{\Dom}{Dom}

\DeclareMathOperator{\supp}{supp}
\newtheorem{theorem}{Theorem}[section]
\newtheorem{lemma}[theorem]{Lemma}
\newtheorem{proposition}[theorem]{Proposition}

\newtheorem{fact}[theorem]{Fact}
\newtheorem{question}[theorem]{Question}
\theoremstyle{definition}
\newtheorem{definition}[theorem]{Definition}
\theoremstyle{remark}

\numberwithin{equation}{section}
\setlist[enumerate]{label=\textup{(\roman*)},leftmargin=*,itemsep=3pt}
\title[Finiteness of  $L^2$ Dolbeault cohomology]{The finiteness of $L^2$ Dolbeault cohomology groups}
\author{Jian Chen}
\address{School of Mathematics and Statistics, Central China Normal University,
Wuhan 430079, People's Republic of China}
\email{jian-chen@whu.edu.cn}
\thanks{The first named author is supported by  NSFC  (Grant No. 12601144).   The second named author is supported by Hainan University (Grant No. XJ2600003643).}
\author{Yi Li}
\address{School of Mathematics and Statistics, Hainan University, Haikou, Hainan 570228, People's Republic of China}
\email{ustcyili@ustc.edu.cn}

\date{\today}
\begin{document}  
\begin{abstract} 
In this paper, motivated by J.-P. Demailly's finiteness argument, we establish a criterion for a closed linear operator between Hilbert spaces to have closed range of finite codimension. We then apply this criterion to give affirmative answers to two questions of T. Ohsawa on the finiteness of $L^2$ Dolbeault cohomology groups in a unified manner.
To this end, we adapt Demailly's compactness argument to a family of complete approximating metrics, thereby establishing the finite dimensionality of the $L^2$ Dolbeault cohomology groups  with respect to a possibly incomplete  metric.  Furthermore, we extend the answers to the setting of singular  metrics with quasi-plurisubharmonic local weights.

\end{abstract}
\maketitle

\setcounter{tocdepth}{1}
\tableofcontents

\section{Introduction}

As is well-known, for a compact complex manifold, elliptic theory ensures
finite-dimensional Dolbeault cohomology without any
curvature assumption, while positivity of the coefficient
bundle yields various vanishing theorems in the compact
K\"ahler setting.
On noncompact manifolds, ordinary and $L^2$ Dolbeault
cohomology can behave quite differently, since the latter
imposes a global square-integrability condition depending
on the chosen metrics. Consequently, finiteness (or vanishing) theorems for (ordinary or $L^2$) Dolbeault  cohomology on noncompact manifolds under positivity or convexity assumptions are long-standing central topics in complex analysis and geometry (e.g., \cite{AG62},  \cite{Kohn63,Kohn64}).

Among many such works,  the works most closely related to the present paper are those of Ohsawa \cite{Oh79,Oh81,Oh82}, Abdelkader \cite{Ab85}, and Demailly \cite{Dem87}.
More precisely, 
Ohsawa \cite{Oh79,Oh81} proved the finite dimensionality of the ordinary Dolbeault cohomology groups $H^{p,q}(X,E)$ for $p+q>n$ on weakly $1$-complete manifolds when the holomorphic line bundle $E$ is positive outside a compact set. Ohsawa \cite{Oh82} and Abdelkader    
\cite{Ab85} also established finiteness results in general bidegrees under semipositivity and curvature rank conditions outside a compact set, assuming that $X$ is K\"ahler there. Demailly \cite{Dem87} subsequently gave a simplified proof of these finiteness results.

\begin{theorem}[{\cite{Oh79,Oh81,Oh82,Ab85}}]
\label{lem-weak-finiteness}
Let $X$ be a weakly $1$-complete complex manifold of dimension
$n$, and let $(E,h)$ be a smooth Hermitian holomorphic line bundle
on $X$.
Suppose that $X$ admits a smooth Hermitian metric $\alpha$ and
a compact subset $Y\subset X$ such that $\alpha$ is K\"ahler
on $X\setminus Y$ and $\Theta_h$ is semipositive of rank at
least $s$ there, where $1\leq s\leq n$.
Then $\dim_{\mathbb{C}}H^{p,q}(X,E)<\infty$ for
$0\leq p,q\leq n$ satisfying $p+q\geq 2n-s+1$.
\end{theorem}

Note that the  weak $1$-completeness assumption in Theorem \ref{lem-weak-finiteness} cannot  be omitted: let $X=\mathbb{C}^2\setminus\{0\}$, and equip the trivial
line bundle $E$ with the Hermitian metric
$|1|_h^2=e^{-|z|^2}$. Its curvature form
is a K\"ahler form, yet
$\dim_{\mathbb{C}}H^{2,1}(X,E)=\infty$ (for comparison, this example satisfies the assumptions of
Theorem \ref{thm-main-intro} with $\omega_g=\Theta_h$, and thus the corresponding $L^2$ cohomology group is finite dimensional).

In connection with Theorem \ref{lem-weak-finiteness},  Ohsawa \cite[Question 1]{Oh22} raised the following  
question (which is formulated in the notation
of this paper).

\begin{question}[{Ohsawa's question \cite[Question 1]{Oh22}}]\label{ques-ohsawa-intro}
Let $M$ be a connected complex manifold of dimension
$n\geq 1$, and let $(E,h)$ be a smooth Hermitian
holomorphic line bundle on $M$ with the curvature
form denoted by $\Theta_h$. Let $g$ be a smooth Hermitian metric  with the fundamental form denoted by 
$\omega_g$. 
Assume that $M$ admits a complete K\"ahler metric $\widehat\omega$ and
that
\[
 \supp(\omega_g-\Theta_h)\Subset M.
\]
Is the $L^2$ Dolbeault cohomology group (see \S \ref{fa-preliminaries} for notation and conventions)
\[
 H_{(2),{\rm max}}^{n,q}(M,E;g,h)
 :=
 \ker(\bar\partial_g^{n,q})/\operatorname{Im}(\bar\partial_g^{n,q-1})
\]
finite dimensional for each $1\leq q\leq n$? 
\end{question}

Concerning Question \ref{ques-ohsawa-intro}, Ohsawa
\cite[Theorem 2.4]{Oh26} recently established a related finiteness
result under another type of assumptions. In particular, if the
Hermitian metric $g$ in Question \ref{ques-ohsawa-intro} is further
assumed to be complete, then \cite[Theorem 2.4]{Oh26} gives an
affirmative answer to Question \ref{ques-ohsawa-intro}.

Ohsawa describes Question \ref{ques-ohsawa-intro} as a difficult problem, and he then raises several other questions, expected to be comparatively easier, under different geometric assumptions (\cite[Questions 2--6]{Oh22}). One of them is the following question for weakly $1$-complete manifolds,  which may be viewed as a natural $L^2$ analogue of Theorem \ref{lem-weak-finiteness}.

\begin{question}[{Ohsawa's question \cite[Question~2]{Oh22}}]\label{ques-ohsawa-intro-2}
Let $M$ be a connected weakly $1$-complete (not necessarily K\"ahler) complex manifold of dimension
$n\geq 1$, i.e., $M$ admits a smooth plurisubharmonic exhaustion function.
Let $(E,h)$ be a smooth Hermitian holomorphic line bundle on $M$ with
curvature form $\Theta_h$, and let $g$ be a smooth Hermitian metric on $M$
with fundamental form $\omega_g$. Assume that
$\supp(\omega_g-\Theta_h)\Subset M$.
Is the  $L^2$ Dolbeault cohomology group
$H_{(2),\mathrm{max}}^{n,q}(M,E;g,h)$
finite dimensional for every $1\leq q\leq n$?
\end{question}

In the present paper, the smooth metric cases of Theorem \ref{thm-main} and
Theorem \ref{thm-ohsawa-question-two} give affirmative answers to
Questions \ref{ques-ohsawa-intro} and \ref{ques-ohsawa-intro-2},
respectively.

We utilize a unified method to establish the smooth metric cases of
Theorems \ref{thm-main} and \ref{thm-ohsawa-question-two}.
The key ingredient is
the following criterion for a closed linear operator to have closed range of finite codimension,
formulated in terms of auxiliary operators on varying Hilbert spaces.
Its proof draws on Demailly's finiteness argument
\cite[Lemma 2.2 and the argument following (2.3)]{Dem87}.

\begin{theorem}[{= Theorem \ref{fa-cancellation}}]\label{fa-cancellation-intro}
Let $T_0:\Dom T_0\subset X_0\to Z_0$ be a closed linear operator
between Hilbert spaces. For $0<\varepsilon\le1$, let
$T_\varepsilon:\Dom T_\varepsilon\subset X_\varepsilon
\to Z_\varepsilon$ be linear operators between Hilbert spaces, and
let $I_\varepsilon:Z_0\to Z_\varepsilon$ be linear maps.
Assume that an integer $N\ge0$ has the following properties.
\begin{enumerate}
\item For each $\varepsilon$ there is a linear map
$G_\varepsilon:\C^N\to Z_\varepsilon$
such that for every $f\in Z_0$ one can choose
$u_{\varepsilon,f}\in\Dom T_\varepsilon$ and
$c_{\varepsilon,f}\in\C^N$  with
\begin{equation}\label{intro-finite-error}
 T_\varepsilon u_{\varepsilon,f}
 =I_\varepsilon f+G_\varepsilon c_{\varepsilon,f},
 \qquad
 \sup_{0<\varepsilon\le1}
 \|u_{\varepsilon,f}\|_{X_\varepsilon}<\infty.
\end{equation}

\item\label{intro-cancellation-ii} Whenever $\varepsilon_\nu\downarrow0$, $f_\nu\to f$ in $Z_0$,
and $u_\nu\in\Dom T_{\varepsilon_\nu}$ satisfy
\begin{equation}\label{intro-limit-hypothesis}
 \sup_\nu\|u_\nu\|_{X_{\varepsilon_\nu}}<\infty,
 \qquad
 T_{\varepsilon_\nu}u_\nu=I_{\varepsilon_\nu}f_\nu,
\end{equation}
one has $f\in\operatorname{Im} T_0$.
\end{enumerate}
Then
\[
 \dim_\C(Z_0/\operatorname{Im} T_0)\le N,
 \qquad \operatorname{Im} T_0\text{ is closed in }Z_0.
\]
\end{theorem}

Applying Theorem \ref{fa-cancellation-intro} to the maximal realization of the
$\bar\partial$ operator and approximating the given metric
by complete Hermitian metrics, we first prove the following
Theorem \ref{thm-main-intro} for the case of  smooth  $h$, answering
Question \ref{ques-ohsawa-intro}. A comparison of local $L^2$
complexes then yields the proof of Theorem \ref{thm-main-intro}.

\begin{theorem}[{= Theorem \ref{thm-main} Answer to Question \ref{ques-ohsawa-intro}}]\label{thm-main-intro}  
With the notation and assumptions of
Question \ref{ques-ohsawa-intro}, except that $h$ is a singular
Hermitian metric with quasi-plurisubharmonic local weights
(in this setting, the support condition is understood in the sense
of currents), the $L^2$ cohomology group
$H_{(2),{\rm max}}^{n,q}(M,E;g,h)$ \footnote{Note that, on a noncompact manifold, the cohomology group
$H^q(M,K_M\otimes E\otimes\mathcal I(h))$ is not in general isomorphic to
$H_{(2),{\rm max}}^{n,q}(M,E;g,h)$.} is finite dimensional
for each $1\leq q\leq n$.
\footnote{Note that a comparison of Theorem \ref{thm-main-intro} with the example
following Theorem \ref{lem-weak-finiteness} shows that
ordinary and $L^2$ Dolbeault cohomology can behave quite
differently: the former can be infinite dimensional even
when the latter is finite dimensional.}
In particular, $\operatorname{Im}(\bar\partial_g^{n,q-1})$
is closed in $L^2_{n,q}(M,E;\omega_g,h)$
for each such $q$.
\end{theorem}

Let us describe the strategy of the proof of Theorem \ref{thm-main-intro} more precisely.
For the  smooth $h$ case of Theorem \ref{thm-main-intro}, we first apply Demailly's
Bochner--Kodaira--Nakano identity
\cite[Theorem~2.12]{Dem86}
to obtain a uniform estimate with curvature and torsion
errors supported in a fixed compact set
(see Proposition \ref{geo-basic}).
We then adapt the local compactness step in Demailly's finiteness argument
\cite[Lemma 2.2]{Dem87} to the family of approximating
metrics:
we use uniform local Sobolev estimates and
finite-dimensional approximation in a fixed local
$L^2$ space defined by $g$ to replace this compact error
by finitely many fixed linear functionals, uniformly
in the approximation parameter
(see Proposition \ref{geo-moments}).
Next, we combine these estimates with metric comparison
in top holomorphic degree (Lemma \ref{geo-comparison})
to obtain solutions of the $\bar\partial$-equation modulo
error spaces, with uniform
$L^2$ bounds in the approximating metrics
(see Lemma \ref{geo-augmented}).
Finally, we cancel the errors in suitable linear combinations
(see the proof of Theorem \ref{fa-cancellation})
before passing to local weak limits (Lemma \ref{geo-limit}),
obtaining the desired finite dimensionality.

For the singular $h$ case of Theorem \ref{thm-main-intro}, we compare its local $L^2$ complex with that
of a smooth metric agreeing with $h$ outside a compact set.
The quotient resolves a coherent sheaf with compact support,
which completes the proof of Theorem \ref{thm-main-intro}.
(see Subsection \ref{ss-singular-metrics}).

Using essentially the same method as for Theorem
\ref{thm-main-intro}, we  give an affirmative answer to
Question \ref{ques-ohsawa-intro-2} and then extend it to the
singular metric setting. The only difference in the smooth metric
case is that we use a plurisubharmonic exhaustion of the weakly
$1$-complete manifold to construct the complete approximating
metrics.  The argument for the singular metric case is exactly the same.

\begin{theorem}[{= Theorem \ref{thm-ohsawa-question-two-singular} Answer to Question \ref{ques-ohsawa-intro-2}}]\label{thm-ohsawa-question-two-singular-intro}
With the notation and assumptions of
Question \ref{ques-ohsawa-intro-2}, except that $h$ is a singular
Hermitian metric with quasi-plurisubharmonic local weights
(in this setting, the support condition is understood in the sense
of currents), the $L^2$ cohomology group
$H_{(2),{\rm max}}^{n,q}(M,E;g,h)$ (in Definition \ref{def-ltwo-cohomology}) is finite dimensional
for each $1\leq q\leq n$.
In particular, $\operatorname{Im}(\bar\partial_g^{n,q-1})$
is closed in $L^2_{n,q}(M,E;\omega_g,h)$
for each such $q$.

\end{theorem}

The rest of the paper is organized as follows.
In Section \ref{fa-preliminaries}, we introduce the notation
and  preliminaries.
In Section \ref{fa-limit-section}, we establish an abstract
criterion for a closed linear operator between Hilbert
spaces to have closed range of finite codimension.  
In Section \ref{sec-geometry}, we derive uniform estimates
for a family of complete approximating metrics and apply
this criterion to prove Theorem \ref{thm-main-intro} first for smooth $h$ case.  We then obtain the singular metric case by comparing the corresponding
local and global $L^2$ complexes with those associated with a suitable
smooth Hermitian metric.
In Section \ref{sec-ohsawa-question-two}, we adapt the method developed
in Section \ref{sec-geometry} to the weakly $1$-complete setting and prove Theorem \ref{thm-ohsawa-question-two-singular-intro}.

\section{Preliminaries}\label{fa-preliminaries}

In this section, we present some standard preliminaries and establish the notational convention.
Let $M$ be a connected complex manifold of dimension $n\geq1$, and let
$E$ be a holomorphic line bundle on $M$, equipped with a possibly
singular Hermitian metric $h$ whose local weights are
quasi-plurisubharmonic (i.e., locally the sum of a plurisubharmonic
function and a smooth function) and not identically $-\infty$.
Fix a smooth Hermitian metric $g$ (not necessarily complete) on $M$, denote its fundamental form
by $\omega_g$, and put $dV_g=\omega_g^n/n!$.

\subsection{Curvature conventions}

\begin{definition}\label{def-curvature}
Let $U\subset M$ be a trivialization domain of $E$, and let
$e$ be a holomorphic frame of $E$ over $U$.
Write $|e|_h^2=e^{-\varphi}$, where
$\varphi=-\log|e|_h^2$ is quasi-plurisubharmonic on $U$.
We denote the curvature current by
\[
  \Theta_h:=\sqrt{-1}\partial\bar\partial\varphi,
\]
which is clearly
independent of the local holomorphic frame and thus
 a globally defined real $(1,1)$-current on $M$.
 The multiplier ideal sheaf \emph{$\mathcal I(h)$} is defined to be the
sheaf of holomorphic function germs $f$ such that $|f|^2e^{-\varphi}$ is
locally integrable with respect to the Euclidean  measure
in any local holomorphic coordinates. 
\end{definition}

\subsection{$L^2$ spaces, maximal $\bar\partial$-operators and $L^2$ cohomologies}\label{ss-closed-form target}

\begin{definition}\label{def-ltwo-spaces}
For $0\leq p,q\leq n$, let $L^2_{p,q}(M,E;\omega_g,h)$ be the Hilbert
space of square-integrable $E$-valued $(p,q)$-forms, equipped with the
inner product $(\cdot,\cdot)_g$ induced by $g$ and $h$.
Its norm is given by
\[
  \|u\|_g^2=\int_M |u|_{g,h}^2\,dV_g.
\]
We denote by $\mathcal A_{\rm c}^{p,q}(M,E)$ the space of compactly
supported smooth $E$-valued $(p,q)$-forms.
\end{definition}

Since the local weights of $h$ are locally bounded above, every
form $u$ in $L^2_{p,q}(M,E;\omega_g,h)$ has locally square-integrable
coefficients. Thus  $\bar\partial u$  is well defined
in the sense of distributions.

\begin{definition}\label{maximal-realization}
For $0\leq p\leq n$ and $0\leq q<n$, \emph{the maximal realization/extension of
$\bar\partial$} is the operator (the metric $g$ is not fixed in \S \ref{sec-geometry} and thus we use the subscript)
\[
  \bar\partial_g^{p,q}:
  \operatorname{Dom}(\bar\partial_g^{p,q})
  \subset L^2_{p,q}(M,E;\omega_g,h)
  \longrightarrow L^2_{p,q+1}(M,E;\omega_g,h),
\]
where the domain
\[
  \operatorname{Dom}(\bar\partial_g^{p,q})
  =\bigl\{u\in L^2_{p,q}(M,E;\omega_g,h):
       \bar\partial u\in L^2_{p,q+1}(M,E;\omega_g,h)
       \text{ in the sense of distribution}\bigr\},
\]
and $\bar\partial_g^{p,q}u=\bar\partial u$.
\end{definition}

We now show that each $\bar\partial_g^{p,q}$ is densely
defined and closed (this is well-known if $h$ is smooth). 

\begin{fact}\label{fact-closed-den}
    Each maximal realization $\bar\partial_g^{p,q}$ of $\bar\partial$ is densely
defined and closed.
\end{fact}

\begin{proof}

Denote by $Z$ the zero variety $V(\mathcal I(h))$. Then $Z$ is a
proper analytic subset of $M$, and $e^{-\varphi}$ is locally
integrable on $M\setminus Z$. In particular, since $Z$ has measure
zero, restriction induces a natural isometric isomorphism
\[
 L^2_{p,q}(M,E;\omega_g,h)
 \cong
 L^2_{p,q}(M\setminus Z,E;\omega_g,h).
\]
Moreover, $e^{-\varphi}\in L^1_{\rm loc}(M\setminus Z).$
Thus by a standard cutoff and regularization argument, we can obtain that  
$\mathcal A_{\rm c}^{p,q}(M\setminus Z,E)$ is dense in
$L^2_{p,q}(M\setminus Z,E;\omega_g,h)\cong  L^2_{p,q}(M,E;\omega_g,h)$.

For any
$u\in\mathcal A_{\rm c}^{p,q}(M\setminus Z,E)$,  its extension by zero
is a smooth compactly supported form on $M$. Furthermore,  both $u$ and
$\bar\partial u$ are square integrable with respect to $g$ and
$h$. Therefore
\[
 \mathcal A_{\rm c}^{p,q}(M\setminus Z,E)
 \subseteq\operatorname{Dom}(\bar\partial_g^{p,q}).
\]
Thus
$\bar\partial_g^{p,q}$ is densely defined.

To prove the closedness, let
$u_j\in\operatorname{Dom}(\bar\partial_g^{p,q})$ satisfy
\[
 u_j\longrightarrow u
 \quad\text{and}\quad
 \bar\partial_g^{p,q}u_j\longrightarrow v
\]
in the corresponding weighted $L^2$ spaces. Since the local
weights of $h$ are quasi-psh and thus locally bounded above, these above convergences
imply  convergence in the sense
of distributions. By the continuity of the distributional
$\bar\partial$, we obtain $\bar\partial u=v$. Consequently,
$u\in\operatorname{Dom}(\bar\partial_g^{p,q})$ and
$\bar\partial_g^{p,q}u=v$. This proves that  
$\bar\partial_g^{p,q}$ is closed.
\end{proof}

\begin{definition}\label{def-adjoints}
For $0\leq q<n$, the domain of the \emph{Hilbert adjoint}
$(\bar\partial_g^{n,q})^*$ consists of those
$v\in L^2_{n,q+1}(M,E;\omega_g,h)$ for which there is a $w\in L^2_{n,q}(M,E;\omega_g,h)$ such that
\[
  (\bar\partial_g^{n,q}u,v)_g=(u,w)_g
  \qquad\text{for each }u\in\operatorname{Dom}(\bar\partial_g^{n,q});
\]
in this case $(\bar\partial_g^{n,q})^*v=w$.

When $h$ is smooth, let $\vartheta_{g,q+1}$ denote
\emph{the formal adjoint} from degree $q+1$ to degree $q$,
characterized by
\[
  (\bar\partial u,v)_g=(u,\vartheta_{g,q+1}v)_g
\]
for compactly supported smooth forms $u,v$ of the appropriate degrees.
In this case, $(\bar\partial_g^{n,q})^*v=\vartheta_{g,q+1}v$
in the sense of distributions.
\end{definition}

Fix $1\leq q\leq n$. It is convenient to regard the preceding
differential as an operator with target $\ker(\bar\partial_g^{n,q})$:
\[
  T_g:\operatorname{Dom}(\bar\partial_g^{n,q-1})\subset L^2_{n,q-1}(M,E;\omega_g,h)\longrightarrow \ker(\bar\partial_g^{n,q}),
  \qquad T_gu=\bar\partial_g^{n,q-1}u.
\]
This operator is densely defined and closed. Its Hilbert adjoint is
\begin{equation}\label{fa-target-adjoint}
  \operatorname{Dom}T_g^*=\ker(\bar\partial_g^{n,q})\cap\operatorname{Dom}((\bar\partial_g^{n,q-1})^*),
  \qquad T_g^*v=(\bar\partial_g^{n,q-1})^*v.
\end{equation}
Indeed, for $v\in \ker(\bar\partial_g^{n,q})$, the boundedness condition on
$u\mapsto(T_gu,v)_g$ defining the adjoint is exactly the condition
defining $(\bar\partial_g^{n,q-1})^*$. 

Clearly, the densely defined closed operators $\bar\partial_g^{n,q}$ form a Hilbert space complex.

\begin{definition}\label{def-ltwo-cohomology}
For $1\leq q\leq n$, the \emph{maximal  $L^2$ cohomology groups} are defined by
\[
  \begin{aligned}
    H_{(2),{\rm max}}^{n,q}(M,E;g,h)&:=\ker(\bar\partial_g^{n,q})/\operatorname{Im}(\bar\partial_g^{n,q-1}).
  \end{aligned}
\]
\end{definition}

\begin{definition}\label{def-graph-norm}
Let $T:\operatorname{Dom}T\subset X\to Z$ be a closed linear operator
between Hilbert spaces. Its graph norm is defined by
\[
  \|u\|_{\rm graph}^2=\|u\|_X^2+\|Tu\|_Z^2.
\]  
Clearly, the domain $\operatorname{Dom}T$ is a Hilbert space w.r.t. this norm.
\end{definition}

\subsection{Density in the graph norm}\label{ss-density}

We shall use the following classical density result   w.r.t. the graph norm.

\begin{lemma}[{e.g., \cite[Chapter VIII, Theorem 3.2]{Dem12}}]\label{lem-graph-core}
 Let $\omega$ be a complete smooth
Hermitian metric.  Assume that $h$ is smooth.
Then for $1\leq q\leq n$, the space $\mathcal A_{\rm c}^{n,q}(M,E)$ is dense
in $\operatorname{Dom}(\bar\partial_\omega^{n,q})\cap\operatorname{Dom}((\bar\partial_\omega^{n,q-1})^*)$ w.r.t. the norm
\[
  \|u\|_\omega+\|\bar\partial_\omega^{n,q}u\|_\omega+\|(\bar\partial_\omega^{n,q-1})^*u\|_\omega.
\]
\end{lemma}

\subsection{Exterior multiplication and graded commutators conventions}

Throughout this paper, brackets of homogeneous operators are understood in the graded sense.

\begin{definition}\label{def-graded-commutator}
For a differential form $\alpha$, let $L_\alpha$ denote exterior
multiplication by $\alpha$, so $L_\alpha(\beta)=\alpha\wedge\beta$.
An operator $A$ is \emph{homogeneous of degree $r\in\mathbb Z$} if it maps
forms of degree $k$ to forms of degree $k+r$ for each $k$.
For homogeneous operators $A$ and $B$ of degrees $r$ and $s$,
respectively, their graded commutator is
\[
  [A,B]:=AB-(-1)^{rs}BA.
\]
\end{definition}

If at least one of $r$ and $s$ is even, then $[A,B]=AB-BA$ is the usual
commutator. If both are odd, then $[A,B]=AB+BA$ is the anticommutator.

\subsection{Local $L^2$ complexes}\label{ss-local-complexes}

Let $\mathcal L_h^{n,q}$ be the sheaf of $E$-valued $(n,q)$-forms
$u$ such that $u$ and $\bar\partial u$ are locally square integrable
with respect to $g$ and $h$. It follows from the local $L^2$ existence  theorem  that the corresponding complex is a fine resolution of $K_M\otimes E\otimes\mathcal I(h)$
\[
0\longrightarrow K_M\otimes E\otimes\mathcal I(h)
\longrightarrow\mathcal L_h^{n,0}
\xrightarrow{\bar\partial}\cdots
\xrightarrow{\bar\partial}\mathcal L_h^{n,n}
\longrightarrow0
\]
(e.g., \cite[Proof of Theorem 5.11]{Dem10}).

\section{A criterion for closed range of finite codimension} \label{fa-limit-section}

In this section, we give a criterion under which a closed  operator between Hilbert spaces has closed range of finite codimension. In the next section, we will apply this criterion to the maximal realization of $\bar{\partial}$ to obtain our desired finiteness theorem.

Motivated by Demailly's finiteness argument
\cite[Lemma 2.2, inequality (2.3), and the subsequent argument]{Dem87}
and his treatment of the $\bar\partial$-equation with an error term (\cite[Remark 3.2]{Dem00}),
we  establish the following variant of the standard $L^2$ existence
lemma with a finite-dimensional error term.
An a priori estimate  yields solvability modulo a fixed finite-dimensional
space. Moreover, it provides a norm estimate for both the solution
and the error term.

\begin{lemma}\label{fa-augmented}
Let $X$ and $Z$ be Hilbert spaces, let $T:\Dom T\subset X\to Z$ be a closed densely
defined  linear operator, and let $L:Z\to\C^N$ be a bounded
linear operator (equivalently, a collection of $N$ bounded linear functionals). Suppose that
$q:\Dom T^*\to[0,\infty)$ is a nonnegative quadratic form and that
there exists a constant $C>0$ such that
\begin{equation}\label{fa-augmented-estimate}
 q(v)\le C\bigl(\|T^*v\|_X^2+|Lv|^2\bigr),
 \qquad \forall v\in\Dom T^*.
\end{equation}
Assume that  $f\in Z$ satisfies  that  there exists a constant $D_f\geq 0$ independent of $v$ such that
\begin{equation}\label{fa-datum-estimate}
 |(v,f)_Z|\le D_f q(v)^{1/2},
 \qquad \forall v\in\Dom T^*.
\end{equation}
Then there are
$u\in\Dom T$ and $c\in\C^N$ such that
\begin{equation}\label{fa-augmented-solution}
 f=Tu+L^*c,\qquad
 \|u\|_X^2+|c|^2\le CD_f^2.
\end{equation}
\end{lemma}

\begin{proof}
Define the linear map  $J:\Dom T^*\to X\oplus\C^N$  by
$$Jv=(T^*v,Lv),$$ and define a functional $\Phi$ on
$\operatorname{Im}J$ by $\Phi(Jv)=(v,f)_Z$.   If $Jv=0$, then  we have $$T^* v=0 \text{ and } L v=0.$$
Consequently, it follows from 
\eqref{fa-augmented-estimate} and \eqref{fa-datum-estimate} that 
$(v,f)_Z=0$, and hence that  $\Phi$ is
well defined. 

Again  by \eqref{fa-datum-estimate} and
\eqref{fa-augmented-estimate} we have
\[
 |\Phi(Jv)|
 \le C^{1/2}D_f\|Jv\|_{X\oplus\C^N}.
\]
Thus $\Phi$ is a bounded functional on $\operatorname{Im}J$.
Extend $\Phi$ by the Hahn--Banach theorem to a bounded linear
functional on $X\oplus\C^N$ with the same bound, still denoted by
$\Phi$. By the Riesz representation theorem, there exists
$(u,c)\in X\oplus\C^N$ such that
\[
 \Phi(x,d)=(x,u)_X+(d,c)_{\C^N},
 \qquad
 (x,d)\in X\oplus\C^N,
\]
and
\begin{equation}\label{equa-bound-u-c}
    \|u\|_X^2+|c|^2
 \le CD_f^2.
\end{equation}
By the expression of $\Phi$ on
$\operatorname{Im}J$, we obtain that, for each $v\in\Dom T^*$,  we have
\[
 (v,f)_Z
 =\Phi(Jv)
 =(T^*v,u)_X+(Lv,c)_{\C^N}.
\]
Since $(Lv,c)_{\C^N}=(v,L^*c)_Z$, 
we have
\[
 (v,f)_Z=(T^*v,u)_X+(v,L^*c)_Z,
 \quad v\in\Dom T^*,
\]
and hence $$
\left(T^* v, u\right)_X=\left(v, f-L^* c\right)_Z, \quad v \in \operatorname{Dom} T^*.
$$
Consequently, 
$u\in\Dom(T^*)^*$ and $(T^*)^*u=f-L^*c$.
Since $T$ is closed and densely defined, $(T^*)^*=T$ and hence $f=Tu+L^*c$.
This, together with \eqref{equa-bound-u-c},   completes the proof of Lemma \ref{fa-augmented}.
\end{proof}

Motivated by Demailly's finiteness argument
\cite[Lemma 2.2 and the argument following (2.3)]{Dem87},
we establish Theorem \ref{fa-cancellation} for operators on
varying Hilbert spaces, which plays a central role in proving our finiteness Theorem \ref{thm-main}. In Theorem \ref{fa-cancellation}, 
by cancelling the errors in suitable linear combinations before
passing to the limit, we obtain the codimension bound from
uniform solution bounds and condition \ref{fa-cancellation-ii}.

\begin{theorem}\label{fa-cancellation}
Let $T_0:\Dom T_0\subset X_0\to Z_0$ be a closed linear operator
between Hilbert spaces. For $0<\varepsilon\le1$, let
$T_\varepsilon:\Dom T_\varepsilon\subset X_\varepsilon
\to Z_\varepsilon$ be linear operators between Hilbert spaces, and
let $I_\varepsilon:Z_0\to Z_\varepsilon$ be linear maps.
Assume that there exists an integer $N\ge 0$ such that the following conditions hold:
\begin{enumerate}
\item For each $\varepsilon$ there is a linear map
$G_\varepsilon:\C^N\to Z_\varepsilon$
such that for every $f\in Z_0$ one can choose
$u_{\varepsilon,f}\in\Dom T_\varepsilon$ and
$c_{\varepsilon,f}\in\C^N$  with
\begin{equation}\label{fa-finite-error}
 T_\varepsilon u_{\varepsilon,f}
 =I_\varepsilon f+G_\varepsilon c_{\varepsilon,f},
 \qquad
 \sup_{0<\varepsilon\le1}
 \|u_{\varepsilon,f}\|_{X_\varepsilon}<\infty.
\end{equation}

\item\label{fa-cancellation-ii} Whenever $\varepsilon_\nu\downarrow0$, $f_\nu\to f$ in $Z_0$,
and $u_\nu\in\Dom T_{\varepsilon_\nu}$ satisfy
\begin{equation}\label{fa-exact-limit-hypothesis}
 \sup_\nu\|u_\nu\|_{X_{\varepsilon_\nu}}<\infty,
 \qquad
 T_{\varepsilon_\nu}u_\nu=I_{\varepsilon_\nu}f_\nu,
\end{equation}
one has $f\in\operatorname{Im} T_0$.
\end{enumerate}
Then
\[
 \dim_\C(Z_0/\operatorname{Im} T_0)\le N,
 \qquad \operatorname{Im} T_0\text{ is closed in }Z_0.
\]
\end{theorem}

\begin{proof}
Take $f_0,\ldots,f_N\in Z_0$. For each $i=0,\ldots,N$,
choose $u_{\varepsilon,f_i}$ and $c_{\varepsilon,f_i}$
as in \eqref{fa-finite-error}.
Then  we have
\begin{equation}\label{error-esti-sup-uefi}
T_\varepsilon u_{\varepsilon,f_i}
=I_\varepsilon f_i+G_\varepsilon c_{\varepsilon,f_i},
\qquad
\sup_{0<\varepsilon\le1}
\|u_{\varepsilon,f_i}\|_{X_\varepsilon}<\infty
\end{equation}
for each $i=0,\ldots,N$.

For each $\varepsilon$, since the $N+1$ vectors
$c_{\varepsilon,f_0},\ldots,c_{\varepsilon,f_N}$
lie in $\C^N$, there exists a unit vector
$a_\varepsilon=(a_{0,\varepsilon},\ldots,a_{N,\varepsilon})
\in\C^{N+1}$ such that
\begin{equation}\label{sum=0-ac}
 \sum_{i=0}^N a_{i,\varepsilon}c_{\varepsilon,f_i}=0.   
\end{equation}
Put
\[
 w_\varepsilon=\sum_{i=0}^N a_{i,\varepsilon}u_{\varepsilon,f_i},
 \qquad f_\varepsilon=\sum_{i=0}^N a_{i,\varepsilon}f_i.
\]

By the linearity of $T_\varepsilon$, $I_\varepsilon$, and $G_\varepsilon$, we have
\[
\begin{aligned}
T_\varepsilon w_\varepsilon
&=I_\varepsilon f_\varepsilon
+G_\varepsilon\left(
\sum_{i=0}^N a_{i,\varepsilon}c_{\varepsilon,f_i}
\right).
\end{aligned}
\]
Since the error term vanishes by \eqref{sum=0-ac}, it follows that
\[
T_\varepsilon w_\varepsilon=I_\varepsilon f_\varepsilon.
\]

Moreover, by the  Cauchy--Schwarz inequality, we have
\[
\begin{aligned}
\|w_\varepsilon\|_{X_\varepsilon}
&=\left\|\sum_{i=0}^N a_{i,\varepsilon}u_{\varepsilon,f_i}\right\|_{X_\varepsilon}\\
&\le \sum_{i=0}^N |a_{i,\varepsilon}|
\|u_{\varepsilon,f_i}\|_{X_\varepsilon}\\
&\le \left(\sum_{i=0}^N |a_{i,\varepsilon}|^2\right)^{1/2}
\left(\sum_{i=0}^N
\|u_{\varepsilon,f_i}\|_{X_\varepsilon}^2\right)^{1/2}\\
&\le \left(\sum_{i=0}^N
\sup_{0<\delta\le1}
\|u_{\delta,f_i}\|_{X_\delta}^2\right)^{1/2}<\infty,
\end{aligned}
\]
where the last inequality follows from $|a_\varepsilon|=1$ and  the estimate \eqref{error-esti-sup-uefi}. Since the final bound is independent of $\varepsilon$, $\|w_\varepsilon\|_{X_\varepsilon}$
is bounded uniformly in $\varepsilon$.

Choose any sequence $\varepsilon_\nu\downarrow0$. Since $|a_{\varepsilon_\nu}|=1$ for each $\nu$, the sequence $\{a_{\varepsilon_\nu}\}_\nu$ lies in the unit sphere of $\C^{N+1}$, which is compact. Thus, after passing to a subsequence, we may assume that
\[
a_{\varepsilon_\nu}\longrightarrow a
\quad\text{in }\C^{N+1}
\]
for some $a=(a_0,\ldots,a_N)\in\C^{N+1}$  with  $|a|=1$.
Hence
\[
 f_{\varepsilon_\nu}\longrightarrow\sum_{i=0}^Na_if_i
 \quad\text{in }Z_0.
\]
Consequently, it follows from condition \ref{fa-cancellation-ii}  that 
\[
\sum_{i=0}^N a_i f_i\in\operatorname{Im} T_0.
\]
Thus, 
 we have
\[
\sum_{i=0}^N a_i[f_i]=0
\quad\text{in }Z_0/\operatorname{Im} T_0,
\]
where $[f_i]$ denotes the class of $f_i$ modulo $\operatorname{Im} T_0$.
Since $a\neq 0$,  the classes
$[f_0],\ldots,[f_N]$ are linearly dependent.
As $f_0,\ldots,f_N$ were arbitrary, we conclude that
\begin{equation}\label{finite-dim-leqN}
    \dim_\C(Z_0/\operatorname{Im} T_0)\le N.
\end{equation}

Since $T_0$ is closed, $\operatorname{Dom}T_0$  equipped with the graph norm, is a Banach space, and $T_0$ defines a bounded linear operator from this space to $Z_0$.  Combining the finiteness \eqref{finite-dim-leqN} with  the open mapping theorem, we conclude that $\operatorname{Im}T_0$ is closed in $Z_0$ (e.g., \cite[Exercise 2.27]{Bre11}).  This completes the proof of Theorem \ref{fa-cancellation}. 
\end{proof}

\section{Complete metric approximation  and proof  of Theorem \ref{thm-main-intro}}
\label{sec-geometry}

In this section, we combine the criterion obtained in Section \ref{fa-limit-section}
with complete metric approximation to prove Theorem \ref{thm-main-intro}
for smooth $h$, and then extend it to the singular $h$ case by  comparison with certain smooth metric of the line bundle.

\begin{theorem}[{=Theorem \ref{thm-main-intro}}]\label{thm-main}
Let $M$ be a connected complex manifold of dimension
$n\geq 1$, and let $E$ be a holomorphic line bundle on $M$
with a possibly singular Hermitian metric $h$ whose local weights
are quasi-plurisubharmonic and not identically $-\infty$.
Denote its curvature current by $\Theta_h$.
Let $g$ be a smooth Hermitian metric with the fundamental form 
$\omega_g$. 
Assume that $M$ admits a complete K\"ahler metric $\widehat\omega$ and
that
\[
 \supp(\omega_g-\Theta_h)\Subset M
 \]
in the sense of currents.
Then for each $1\leq q\leq n$, the $L^2$-cohomology group 
\[
 H_{(2),\max}^{n,q}(M,E;g,h)=\ker(\bar\partial_g^{n,q})/\operatorname{Im}(\bar\partial_g^{n,q-1})
\]
is finite dimensional (see \S \ref{fa-preliminaries} for notation and conventions).
In particular, $\operatorname{Im}(\bar\partial_g^{n,q-1})$ is closed in $L^2_{n,q}(M,E;\omega_g,h)$.
\end{theorem}

Note that from now on,  until Subsection \ref{ss-singular-metrics}, all notation and hypotheses
are as in Theorem \ref{thm-main} and Section \ref{fa-preliminaries},
with the additional assumption that $h$ is smooth.

\subsection{Metric comparison in top holomorphic degree}

Note that  $\omega_g$ is K\"ahler outside a compact subset of $M$.
For $0<\varepsilon\leq1$, set
\begin{equation}\label{geo-regularization}
  \omega_\varepsilon:=\omega_g+\varepsilon\widehat\omega.
\end{equation}
Throughout this  Section, subscripts $\varepsilon$ indicate that the corresponding $L^2$
spaces, maximal differentials, and Hilbert adjoints are defined using
$\omega_\varepsilon$ on $M$ and $h$ on $E$.
Thus $L^2_{p,q}(M,E;\omega_\varepsilon,h)$,
$\bar\partial_\varepsilon^{p,q}$, and
$(\bar\partial_\varepsilon^{p,q})^*$ are obtained from
Definitions \ref{def-ltwo-spaces}, \ref{maximal-realization}, and
\ref{def-adjoints} by replacing $\omega_g$ with $\omega_\varepsilon$.
In particular, we write
\[
  \|u\|_\varepsilon^2
  =\int_M |u|_{\omega_\varepsilon,h}^2
     \frac{\omega_\varepsilon^n}{n!}.
\]
Each $\omega_\varepsilon$ is a complete Hermitian metric and hence
Lemma \ref{lem-graph-core} applies to these metrics.

Let $\Lambda_\varepsilon$ be the pointwise adjoint of
$L_{\omega_\varepsilon}$ with respect to the Hermitian inner products
induced by $\omega_\varepsilon$ and $h$.
For $1\leq q\leq n$, we define  the curvature-type endomorphism
\begin{equation}\label{Aepsi-aepsi}
  A_{\varepsilon,q}
  :=[L_{\omega_g},\Lambda_\varepsilon]\big|_{
      \Lambda^{n,q}T^*M\otimes E}.
\end{equation}
The endomorphism $A_{\varepsilon,q}$ is pointwise positive definite,
as verified in the  following Lemma \ref{geo-comparison}.
We denote by $A_{\varepsilon,q}^{-1/2}$ its pointwise positive inverse
square root. For a measurable $E$-valued $(n,q)$-form $f$, we understand
\[
  \|A_{\varepsilon,q}^{-1/2}f\|_\varepsilon^2
  :=\int_M (A_{\varepsilon,q}^{-1}f,f)_{\omega_\varepsilon,h}
       \frac{\omega_\varepsilon^n}{n!},
\]
where the integral is initially allowed to be infinite.
Its finiteness for $f\in L^2_{n,q}(M,E;\omega_g,h)$ will follow from
Lemma \ref{geo-comparison}.

We first give the following observation, which is essentially due to \cite[Proposition 1.1]{Oh84} and \cite[Lemmas 3.2–3.3]{Dem82}.

\begin{lemma}[{\cite{Oh84}, \cite{Dem82}}]\label{geo-comparison}
The operator $A_{\eps,q}$ is pointwise positive and
\begin{equation}\label{Aepsilon-qI}
0<A_{\eps,q}\leq qI.    
\end{equation}
 Furthermore, for every $f\in L^2_{n,q}(M,E;\omega_g,h)$, we have
\begin{equation}\label{geo-datum}
 \|f\|_\eps\leq\|f\|_g, 
 \qquad
 \|A_{\eps,q}^{-1/2}f\|_\eps^2\leq\frac1q\|f\|_g^2.
\end{equation}
 \footnote{The first inequality in \eqref{geo-datum} is obtained in \cite[Proposition 1.1]{Oh84}. It also holds for $q=0$, where it is in fact an equality.}
\end{lemma}

\begin{proof}
Fix $1\leq q\leq n$.
For each $x\in M$, let $\mu_1,\ldots,\mu_n$ be
the eigenvalues of $\omega_g$ with respect to
$\omega_\varepsilon$ at $x$.
Since $\omega_\varepsilon\geq\omega_g>0$, we have
$0<\mu_j\leq1$ for $1\leq j\leq n$.
The operator $A_{\varepsilon,q}$ is self-adjoint with
respect to the Hermitian inner product induced by
$\omega_\varepsilon$ and $h$, and its eigenvalues at $x$
are $\sum_{j\in J}\mu_j$, where $J$ ranges over the
$q$-element subsets of $\{1,\ldots,n\}$.
Each such sum lies in $(0,q]$, and thus $0<A_{\varepsilon,q}\leq qI$.   This proves \eqref{Aepsilon-qI}.

Let $\Lambda_g$ be the adjoint of $L_{\omega_g}$ with
respect to $g$, and set
\[
A_{g,q}
:=
[L_{\omega_g},\Lambda_g]
\big|_{\Lambda^{n,q}T^*M\otimes E}
=
qI.
\]
Write $dV_g:=\omega_g^n/n!$ and
$dV_\varepsilon:=\omega_\varepsilon^n/n!$.

Applying  the metric comparison results in
\cite[Lemmas~3.2--3.3]{Dem82} (or \cite[Chapter VIII, proof of Lemma 6.3, pp. 376--377;
Chapter VI, Corollary 5.9, p. 301]{Dem12}) to
$\omega_g\leq\omega_\varepsilon$,
we obtain the following pointwise inequalities for any
$E$-valued $(n,q)$-form $f$:
\[
\begin{aligned}
|f|_{\omega_\varepsilon,h}^{2}\,dV_\varepsilon
&\leq
|f|_{g,h}^{2}\,dV_g,\\
\big\langle A_{\varepsilon,q}^{-1}f,f
\big\rangle_{\omega_\varepsilon,h}\,dV_\varepsilon
&\leq
\big\langle A_{g,q}^{-1}f,f
\big\rangle_{g,h}\,dV_g
=
\frac1q|f|_{g,h}^{2}\,dV_g.
\end{aligned}
\]
Integrating over $M$ proves \eqref{geo-datum}.
\end{proof}

Recall from   Subsection \ref{ss-closed-form target} that we may regard the  $\bar\partial$ as an operator into the space  of  closed forms:
\begin{equation}\label{geo-target}
 T_\eps:\Dom(\bar\partial_\varepsilon^{n,q-1})\subset L^2_{n,q-1}(M,E;\omega_\varepsilon,h)\longrightarrow \ker(\bar\partial_\varepsilon^{n,q}),
 \qquad T_\eps u=\dbar u.
\end{equation}
It is densely defined and closed, and 
\begin{equation}\label{geo-adjoint}
 \Dom T_\eps^*=\ker(\bar\partial_\varepsilon^{n,q})\cap\Dom((\bar\partial_\varepsilon^{n,q-1})^*),
 \qquad T_\eps^*v=(\bar\partial_\varepsilon^{n,q-1})^*v.
\end{equation}

\subsection{Uniform $L^2$ estimates with a cutoff term}

As shown in Lemma \ref{geo-comparison},
$A_{\eps,q}$ defines a bounded positive
operator on $L^2_{n,q}(M,E;\omega_\varepsilon,h)$ for $q\geq1$.
Thus $v\mapsto(A_{\eps,q}v,v)_\eps$ is a nonnegative
continuous quadratic form on $L^2_{n,q}(M,E;\omega_\varepsilon,h)$.

We now utilize curvature positivity outside a compact set
through the Hermitian Bochner--Kodaira--Nakano identity, while
controlling the possible error on the compact set by a cutoff term.
This leads to the following uniform $L^2$ estimate.

\begin{proposition}\label{geo-basic}
There exist $\chi\in C_c^\infty(M)$ and $C\geq0$, independent of
$\eps\in(0,1]$, such that
\begin{equation}\label{geo-basic-estimate}
 (A_{\eps,q}v,v)_\eps\leq\|\bar\partial_\varepsilon^{n,q}v\|_\eps^2
       +\|(\bar\partial_\varepsilon^{n,q-1})^*v\|_\eps^2+C\|\chi v\|_\eps^2
\end{equation}
for all $E$-valued $(n,q)$-forms
$v\in\Dom(\bar\partial_\varepsilon^{n,q})\cap\Dom((\bar\partial_\varepsilon^{n,q-1})^*)$.
In particular, 
\begin{equation}\label{geo-closed-estimate}
(A_{\eps,q}v,v)_\eps\leq
\|T_\eps^*v\|_\eps^2+C\|\chi v\|_\eps^2
\end{equation}
for any $v\in \Dom T_\eps^*$.
\end{proposition}

\begin{proof}

Let $B=\Theta_h-\omega_g$. By hypothesis, $B$ has compact support.
On compactly supported smooth bundle-valued forms, Demailly's
 Bochner--Kodaira--Nakano identity  \cite[Theorem~2.12]{Dem86}
(or see \cite[Chapter VII, Theorems 1.1, 1.4]{Dem12}) gives
\begin{equation}\label{geo-bkn}
 \Delta_\varepsilon''
 =
 \Delta_{\varepsilon,\tau}'
 +[L_{\Theta_h},\Lambda_\varepsilon]
 +\mathcal T_{\omega_\varepsilon},
\end{equation}
where
\[
 \Delta_{\eps,\tau}'=
 [D_h'+\tau_\eps,(D_h'+\tau_\eps)^{*,\eps}],
 \qquad
 \tau_\eps=[\Lambda_\eps,L_{\partial\omega_\eps}],
\]
\begin{equation}\label{geo-torsion}
 \mathcal T_{\omega_\eps}=
 \left[\Lambda_\eps,
   \left[\Lambda_\eps,L_{\frac{\sqrt{-1}}2
                           \partial\bar\partial\omega_\eps}\right]\right]
 -[L_{\partial\omega_\eps},L_{\partial\omega_\eps}^{*,\eps}],
\end{equation}
and $D_h'$ denotes the $(1,0)$ part of the Chern connection.

\step{1} 
Estimate     for compactly supported
smooth $E$-valued $(n,q)$-forms.

We next show that the curvature error
$[L_B,\Lambda_\eps]$ and the torsion term
$\mathcal T_{\omega_\eps}$ can be controlled uniformly
by a compactly supported $L^2$ term.

Since $\widehat{\omega}$ and $\Theta_h$ are $d$-closed, the
definition of $\omega_\varepsilon$ and the equality
$B=\Theta_h-\omega_g$ give
\[
 \partial\omega_\varepsilon
 =
 \partial\omega_g
 =
 -\partial B,
 \qquad
 \partial\bar\partial\omega_\varepsilon
 =
 \partial\bar\partial\omega_g
 =
 -\partial\bar\partial B.
\]

Choose a compact set $K\subset M$ that contains $\operatorname{supp}B$ (independent of $\varepsilon$).
Since $B=0$ on $M\setminus K$, we have
$\omega_g=\Theta_h$ on $M\setminus K$. Since both
$\Theta_h$ and $\widehat\omega$ are $d$-closed, 
$\omega_\varepsilon=\Theta_h+\varepsilon\widehat\omega$
is K\"ahler on $M\setminus K$. 
In particular,
\[
\partial\omega_\varepsilon
=\partial\bar\partial\omega_\varepsilon=0
\qquad\text{on }M\setminus K.
\]
It then follows from \eqref{geo-torsion} that
$\mathcal T_{\omega_\varepsilon}=0$ on $M\setminus K$.
Moreover, $B=0$  implies
$[L_B,\Lambda_\varepsilon]=0$.

Choose a relatively compact open neighborhood $U$ of $K$  and set $\omega_0:=\omega_g$.
For $0\leq\varepsilon\leq1$, the coefficient matrices of
$\omega_\varepsilon$ and their inverses are uniformly bounded
on $\overline U$. The derivatives in \eqref{geo-torsion} are also uniformly bounded on $\overline U$.
Since $B$ and $h$ are fixed and smooth, the pointwise operator
norms of the zero-order self-adjoint endomorphisms
\[
 [L_B,\Lambda_\varepsilon]
 +
 \mathcal T_{\omega_\varepsilon}
\]
are bounded by a constant $C$ independent of $\varepsilon$ on $\overline U$ for $0\leq\varepsilon\leq1$.

Choose $\chi\in C_c^\infty(U,\mathbb R)$ such that
$0\leq\chi\leq1$ and $\chi=1$ on a neighborhood of $K$.
The support property and the uniform bound just established give
\begin{equation}\label{geo-T-error}
    \begin{aligned}
 &\big(
   ([L_B,\Lambda_\varepsilon]
      +\mathcal T_{\omega_\varepsilon})v,v
  \big)_\varepsilon\\
 \geq&
 -C\int_K |v|_{\omega_\varepsilon,h}^2\,dV_{\omega_\varepsilon}
 \geq
 -C\|\chi v\|_\varepsilon^2.
 \end{aligned}
\end{equation}

Apply \eqref{geo-bkn} to a compactly supported smooth
$E$-valued $(n,q)$-form $v$.
Taking the $L^2$ inner product with $v$ and integrating by parts,
we obtain
\begin{equation}\label{geo-bkn-applica}
   \begin{aligned}
 &\|\bar\partial_\varepsilon^{n,q}v\|_\varepsilon^2
  +\|(\bar\partial_\varepsilon^{n,q-1})^*v\|_\varepsilon^2\\
 =&
 (\Delta_\varepsilon''v,v)_\varepsilon\\
 =&
 (\Delta_{\varepsilon,\tau}'v,v)_\varepsilon
 +\big([L_{\Theta_h},\Lambda_\varepsilon]v,v\big)_\varepsilon
 +(\mathcal T_{\omega_\varepsilon}v,v)_\varepsilon.
 \end{aligned} 
\end{equation}

Clearly, any
compactly supported smooth bundle-valued form $v$ satisfies
\[
(\Delta_{\varepsilon,\tau}'v,v)_\varepsilon
=\|(D_h'+\tau_\varepsilon)v\|_\varepsilon^2
+\|(D_h'+\tau_\varepsilon)^{*,\varepsilon}v\|_\varepsilon^2\geq0.
\]
By the definitions of $B$ and $A_{\varepsilon,q}$,  we have
\[
 [L_{\Theta_h},\Lambda_\varepsilon]
 =
 A_{\varepsilon,q}+[L_B,\Lambda_\varepsilon].
\]
Substituting this decomposition into \eqref{geo-bkn-applica},
discarding the nonnegative term
$(\Delta_{\varepsilon,\tau}'v,v)_\varepsilon$, and using the
estimate \eqref{geo-T-error} for the curvature error and the torsion term, we obtain
\[
 \begin{aligned}
 &\|\bar\partial_\varepsilon^{n,q}v\|_\varepsilon^2
  +\|(\bar\partial_\varepsilon^{n,q-1})^*v\|_\varepsilon^2\\
 \geq&
(A_{\eps,q}v,v)_\eps
 +\big(
   ([L_B,\Lambda_\varepsilon]
      +\mathcal T_{\omega_\varepsilon})v,v
  \big)_\varepsilon\\
 \geq&
(A_{\eps,q}v,v)_\eps -C\|\chi v\|_\varepsilon^2.
 \end{aligned}
\]
This proves \eqref{geo-basic-estimate} for compactly supported
smooth $E$-valued $(n,q)$-forms.

\step{2}  Estimate for elements in $\Dom(\bar\partial_\varepsilon^{n,q})\cap\Dom((\bar\partial_\varepsilon^{n,q-1})^*)$.

  To pass the estimate for compactly supported
smooth forms to general elements in $\Dom(\bar\partial_\varepsilon^{n,q})\cap\Dom((\bar\partial_\varepsilon^{n,q-1})^*)$, we fix $\varepsilon\in(0,1]$ and let
\[
 v\in
 \operatorname{Dom}(\bar\partial_\varepsilon^{n,q})
 \cap
 \operatorname{Dom}((\bar\partial_\varepsilon^{n,q-1})^*).
\]
Since $\omega_\varepsilon$ is complete,  Lemma \ref{lem-graph-core}   provides a sequence of compactly supported smooth
forms $v_j$, of the same bidegree and with values in the same
bundle as $v$, such that
\[
 \begin{aligned}
 v_j
 &\longrightarrow v,\\
 \bar\partial_\varepsilon^{n,q}v_j
 &\longrightarrow \bar\partial_\varepsilon^{n,q}v,\\
 (\bar\partial_\varepsilon^{n,q-1})^*v_j
 &\longrightarrow (\bar\partial_\varepsilon^{n,q-1})^*v
 \end{aligned}
\]
in the corresponding $L^2$ spaces defined by
$\omega_\varepsilon$ and $h$.

The estimate
$0\leq A_{\varepsilon,q}\leq q\operatorname{I}$ in \eqref{Aepsilon-qI} of  Lemma \ref{geo-comparison}
implies
\[
 \begin{aligned}
\big|(A_{\varepsilon,q}v_j,v_j)_\varepsilon
     -(A_{\varepsilon,q}v,v)_\varepsilon\big|
 \leq
 q\big(\|v_j\|_\varepsilon+\|v\|_\varepsilon\big)
   \|v_j-v\|_\varepsilon
 \longrightarrow0.
 \end{aligned}
\]
Moreover, since $0\leq\chi\leq1$, we have
$\|\chi(v_j-v)\|_\varepsilon
 \leq\|v_j-v\|_\varepsilon\to0$.
Thus all terms in \eqref{geo-basic-estimate}, applied to $v_j$,
converge to the corresponding terms for $v$.
Passing to the limit proves \eqref{geo-basic-estimate} for all
\[
 v\in
 \operatorname{Dom}(\bar\partial_\varepsilon^{n,q})
 \cap
 \operatorname{Dom}((\bar\partial_\varepsilon^{n,q-1})^*).
\]
This completes the proof of  Proposition \ref{geo-basic}.

\end{proof}

\subsection{Uniform estimates with finitely many fixed linear functionals}

Note that Proposition \ref{geo-basic} controls the quadratic form $(A_{\eps,q}v,v)_\eps$ via $\|T_\eps^*v\|_\eps^2$ up to a compactly supported $L^2$ error term. In Proposition \ref{geo-moments} below, we will replace  this error with finitely many linear functionals.

The proof of Proposition \ref{geo-moments} mainly adapts the local
compactness step in Demailly's finiteness argument
(\cite[Lemma 2.2]{Dem87}) to the family of metrics
$\omega_\eps=\omega_g+\eps\widehat\omega$.
A key difference  is that our metrics
vary even on the compact set supporting the error, while the
metric in \cite[Lemma 2.2]{Dem87} is kept fixed on the
corresponding region.
We then establish uniform local Sobolev estimates and utilize
the Rellich compactness argument in a fixed local $L^2$ space
defined by $g$.
This yields control of $(A_{\varepsilon,q}v,v)_\varepsilon$
without requiring a global positive lower bound for
$A_{\varepsilon,q}$.

\begin{proposition}\label{geo-moments}
There exist an integer $N\geq0$, smooth compactly supported
$E$-valued $(n,q)$-forms $\psi_1,\ldots,\psi_N$,
and constants $C_q,C_q'>0$ such that, setting
\[
\ell_j(v):=(v,\psi_j)_g,\qquad j=1,\ldots,N,
\]
the following estimates hold for each $\eps\in(0,1]$:
\begin{equation}\label{geo-moment-estimate}
(A_{\eps,q}v,v)_\eps\leq C_q\left(
\|T_\eps^*v\|_\eps^2+
\sum_{j=1}^N|\ell_j(v)|^2
\right)
\end{equation}
for any $v\in\Dom T_\eps^*$, and
\begin{equation}\label{fucn-esti-j1l}
    \sum_{j=1}^N|\ell_j(u)|^2
\leq C_q'\|u\|_\eps^2
\end{equation}
for any $u\in L^2_{n,q}(M,E;\omega_\varepsilon,h)$.
\end{proposition}

\begin{proof}

\step{1} Uniform local Sobolev estimate.

If $M$ is noncompact, choose a smoothly bounded open set
$U\Subset M$ containing $\operatorname{supp}\chi$.
Choose $\chi_0,\chi_1\in C_c^\infty(U,\mathbb R)$ such that
$0\leq\chi_0,\chi_1\leq1$, $\chi_0=1$ near
$\operatorname{supp}\chi$, and $\chi_1=1$ near
$\operatorname{supp}\chi_0$.
If $M$ is compact, take $U=M$ and $\chi_0=\chi_1=1$.

Choose $C_U\geq0$ such that
$\widehat\omega\leq C_U\omega_g$ on $\overline U$.
It follows from the definition of $\omega_\varepsilon$ and the calculation
in the proof of Lemma \ref{geo-comparison} that 
\begin{equation}\label{geo-local-uniform}
 \omega_g\leq\omega_\varepsilon\leq(1+C_U)\omega_g,
 \qquad
 A_{\varepsilon,q}\geq\frac{q}{1+C_U}\operatorname{Id}
\end{equation}
on $\overline U$ for $0\leq\varepsilon\leq1$ ($\omega_0:=\omega_g$).
In particular, the local $L^2$ norms are uniformly equivalent.

Put $F:=\Lambda^{n,q}T^*M\otimes E$.
We fix a smooth Hermitian connection $\nabla$ on $F$ and use the following Sobolev norm
\[
 \|w\|_{H_g^1(U)}^2
 :=\|w\|_{L_g^2(U)}^2+\|\nabla w\|_{L_g^2(U)}^2,
\]
where all norms in this definition are w.r.t. $g$ and $h$.
Let $H_0^1(U,F)$ denote the closure of $C_c^\infty(U,F)$
with respect to this norm (when $U=M$ is compact, then
$H_0^1(M,F)=H^1(M,F)$).

Let $\vartheta_{\varepsilon,q}$ denote the formal adjoint of
$\bar\partial$ with respect to $\omega_\varepsilon$ and $h$
(from degree $q$ to degree $q-1$). Note that the first-order system
$P_\varepsilon=(\bar\partial,\vartheta_{\varepsilon,q})$
has injective principal symbol.
Indeed, for a nonzero real covector $\xi$ at $x\in\overline U$
and $w\in F_x$, we have
\[
 \begin{aligned}
 &|\xi^{0,1}\wedge w|_{\omega_\varepsilon,h}^2
 +|(\xi^{0,1}\wedge\bullet)^{*,\varepsilon}w|_{\omega_\varepsilon,h}^2\\
 =&
 |\xi^{0,1}|_{\omega_\varepsilon}^2|w|_{\omega_\varepsilon,h}^2\\
 \geq& c|\xi|_g^2|w|_{\omega_\varepsilon,h}^2,
 \end{aligned}
\]
where $c>0$ is independent of $x$ and $\varepsilon$.
Since the coefficients of $P_\varepsilon$ depend smoothly on
$(x,\varepsilon)\in\overline U\times[0,1]$,  the bounds of each 
 coefficient  needed for the interior elliptic estimate
are uniform in $\varepsilon$. Combining this with the uniform symbol bound, we then obtain the following interior elliptic estimate
\begin{equation}\label{geo-local-garding}
 \|\chi_0v\|_{H_g^1(U)}^2
 \leq C\left(
 \|\chi_1\bar\partial v\|_\varepsilon^2
 +\|\chi_1\vartheta_{\varepsilon,q}v\|_\varepsilon^2
 +\|\chi_1v\|_\varepsilon^2\right)
\end{equation}
for smooth $E$-valued $(n,q)$-forms $v$, where $C$ is independent of $\varepsilon$.

Fix $\varepsilon\in(0,1]$ and
$v\in\operatorname{Dom}T_\varepsilon^*$.
By \eqref{geo-adjoint},
\[
 \bar\partial v=0,
 \qquad
 \vartheta_{\varepsilon,q}v=T_\varepsilon^*v
\]
in the sense of distributions.
Since $\omega_\varepsilon$ is complete, by Lemma
\ref{lem-graph-core}, we can choose  forms $v_\nu\in C_c^\infty(M,F)$ such that
\[
 \begin{aligned}
 v_\nu&\longrightarrow v,\\
 \bar\partial v_\nu&\longrightarrow0,\\
 \vartheta_{\varepsilon,q}v_\nu&\longrightarrow T_\varepsilon^*v
 \end{aligned}
\]
in the corresponding $L^2$ spaces defined by
$\omega_\varepsilon$ and $h$.
Applying \eqref{geo-local-garding} to $v_\nu-v_\mu$, we obtain  that
$\chi_0v_\nu$ is a Cauchy sequence in $H_g^1(U)$.
Moreover, $\chi_0v_\nu$ converges to $\chi_0v$ in $L_g^2(U)$.
Thus $\chi_0v\in H_0^1(U,F)$, and via passing to the limit,
we obtain \eqref{geo-local-garding} for
$v\in\operatorname{Dom}T_\varepsilon^*$.
Moreover, \eqref{geo-local-uniform} gives
\[
 \|\chi_1v\|_\varepsilon^2
 \leq\frac{1+C_U}{q}(A_{\eps,q}v,v)_\eps.
\]
We then obtain
\begin{equation}\label{geo-sobolev}
 \|\chi_0v\|_{H_g^1(U)}^2
 \leq C_1\left(\|T_\varepsilon^*v\|_\varepsilon^2
                  +(A_{\eps,q}v,v)_\eps\right),
\end{equation}
where $C_1>0$ is independent of $\varepsilon$.

\step{2} Construction of fixed linear functionals.

Since $\chi_0=1$ on $\operatorname{supp}\chi$, by
\eqref{geo-closed-estimate} and the uniform equivalence of the
local $L^2$ norms following \eqref{geo-local-uniform}, we obtain
\begin{equation}\label{geo-local-error}
(A_{\eps,q}v,v)_\eps
 \leq\|T_\varepsilon^*v\|_\varepsilon^2
      +C_0\|\chi_0v\|_{L_g^2(U)}^2,
\end{equation}
where $C_0\geq0$ is independent of $\varepsilon$.

Choose a sufficiently small $\delta>0$ such that $C_0C_1\delta\leq1/2$.
By the Rellich compactness theorem (e.g., \cite[Theorem A.3.1]{MM07}), the unit ball of
$H_0^1(U,F)$ is relatively compact (and thus  totally bounded) in $L_g^2(U,F)$.
Since $C_c^\infty(U,F)$ is dense in $L_g^2(U,F)$,
we can choose finitely many sections
$f_1,\ldots,f_k\in C_c^\infty(U,F)$ such that for each
$w\in H_0^1(U,F)$ with $\|w\|_{H_g^1(U)}\leq1$,
there exists $1\leq j\leq k$ satisfying
\[
 \|w-f_j\|_{L_g^2(U)}<\sqrt{\delta}.
\]

Let $V$ be the subspace spanned by $f_1,\ldots,f_k$, and let
$P_V$ be the orthogonal projection from $L_g^2(U,F)$ onto $V$.
By the minimizing property of the orthogonal projection, we obtain
\[
 \|w-P_Vw\|_{L_g^2(U)}^2
 \leq\delta\|w\|_{H_g^1(U)}^2
\]
for each $w\in H_0^1(U,F)$.

Choose an $L_g^2(U,F)$-orthonormal basis
$e_1,\ldots,e_N$ of $V$.
By the orthogonal decomposition
$w=P_Vw+(w-P_Vw)$, we obtain
\begin{equation}\label{geo-spectral-tail}
 \|w\|_{L_g^2(U)}^2
 \leq\sum_{j=1}^N|(w,e_j)_{L_g^2(U)}|^2
       +\delta\|w\|_{H_g^1(U)}^2.
\end{equation}

For $1\leq j\leq N$, put $$\psi_j:=\chi_0e_j$$ on $U$.
Extend $\psi_j$ by zero to $M$.
Since $V\subset C_c^\infty(U,F)$,
we have $\psi_j\in C_c^\infty(M,F)$.
By the definition of $\ell_j$, we obtain
\[
 \ell_j(v)=(v,\psi_j)_g=(\chi_0v,e_j)_{L_g^2(U)}.
\]

\step{3} Proof of the two estimates.

Applying \eqref{geo-spectral-tail} to $w=\chi_0v$
and using \eqref{geo-local-error} and \eqref{geo-sobolev},
we obtain
\[
 \begin{aligned}
 (A_{\eps,q}v,v)_\eps
 &\leq\|T_\varepsilon^*v\|_\varepsilon^2
      +C_0\sum_{j=1}^N|\ell_j(v)|^2
      +C_0\delta\|\chi_0v\|_{H_g^1(U)}^2\\
 &\leq\frac32\|T_\varepsilon^*v\|_\varepsilon^2
      +C_0\sum_{j=1}^N|\ell_j(v)|^2
      +\frac12(A_{\eps,q}v,v)_\eps.
 \end{aligned}
\]
Subtracting $\frac12(A_{\eps,q}v,v)_\eps$ from both sides
of the above inequality and multiplying by $2$, we obtain the desired
\eqref{geo-moment-estimate} with $C_q=\max\{3,2C_0\}$.

Finally, for $u\in L^2_{n,q}(M,E;\omega_\varepsilon,h)$, by Bessel's inequality and the uniform equivalence of the
local $L^2$ norms induced by \eqref{geo-local-uniform}, we obtain
\[
 \sum_{j=1}^N|\ell_j(u)|^2
 \leq\|\chi_0u\|_{L_g^2(U)}^2
 \leq C_2\|\chi_0u\|_\varepsilon^2
 \leq C_2\|u\|_\varepsilon^2,
\]
where $C_2>0$ is independent of $\varepsilon$.
Taking $C_q'=C_2$, we obtain the second desired estimate \eqref{fucn-esti-j1l}.  This completes the proof of Proposition \ref{geo-moments}.
\end{proof}

\subsection{Solving $\bar\partial$ with a finite-dimensional error term}

We now use the 
functionals constructed in Proposition \ref{geo-moments} to construct
 a bounded  operator $L_\varepsilon$ and use its 
Hilbert adjoint $L_\varepsilon^*$ to construct a
finite-dimensional error term.
For each $\varepsilon\in(0,1]$,  we restrict the functionals
$\ell_1,\ldots,\ell_N$ in Proposition \ref{geo-moments}
to $\ker(\bar\partial_\varepsilon^{n,q})$ and define
\[
L_\varepsilon:\ker(\bar\partial_\varepsilon^{n,q})\longrightarrow\mathbb C^N,
 \qquad
 L_\varepsilon v=(\ell_1(v),\ldots,\ell_N(v)).
\]
By Proposition \ref{geo-moments},  we obtain that $L_\varepsilon$ is bounded
with the operator norm at most $\sqrt{C_q'}$.
The adjoint
\[
 L_\varepsilon^*:
 \mathbb C^N\longrightarrow \ker(\bar\partial_\varepsilon^{n,q})
\]
is taken with respect to $(\cdot,\cdot)_\varepsilon$ on
$\ker(\bar\partial_\varepsilon^{n,q})$ and the standard Hermitian inner product
on $\mathbb C^N$.
It is bounded with the operator norm at most $\sqrt{C_q'}$, and
$\dim_{\mathbb C}\operatorname{Im}L_\varepsilon^*\leq N$.

\begin{lemma}\label{geo-augmented}
For each $f\in \ker(\bar\partial_g^{n,q})$ and $\varepsilon\in(0,1]$, there exist
$u_{\varepsilon,f}\in\operatorname{Dom}T_\varepsilon$
and $c_{\varepsilon,f}\in\mathbb C^N$ such that
\begin{equation}\label{geo-augmented-equation}
 \begin{aligned}
 &\quad f=T_\varepsilon u_{\varepsilon,f}
       +L_\varepsilon^*c_{\varepsilon,f},\\
 &\|u_{\varepsilon,f}\|_\varepsilon^2
       +|c_{\varepsilon,f}|^2
 \leq C_q\|f\|_g^2.
 \end{aligned}
\end{equation}
\end{lemma}

\begin{proof}
By Lemma \ref{geo-comparison}, $f\in L^2_{n,q}(M,E;\omega_\varepsilon,h)$.
For $v\in\operatorname{Dom}T_\varepsilon^*$,
the Cauchy--Schwarz inequality and \eqref{geo-datum} give the estimate
\[
 \begin{aligned}
 |(v,f)_\varepsilon|^2
 &\leq (A_{\eps,q}v,v)_\eps
       \|A_{\varepsilon,q}^{-1/2}f\|_\varepsilon^2\\
 &\leq\frac1q\|f\|_g^2(A_{\eps,q}v,v)_\eps.
 \end{aligned}
\]
Combining this with \eqref{geo-moment-estimate}, we obtain
\[
 |(v,f)_\varepsilon|^2
 \leq\frac{C_q}{q}\|f\|_g^2
 \left(\|T_\varepsilon^*v\|_\varepsilon^2
       +|L_\varepsilon v|^2\right).
\]
Lemma \ref{fa-augmented} gives
$f=T_\varepsilon u_{\varepsilon,f}
   +L_\varepsilon^*c_{\varepsilon,f}$ with
\[
 \|u_{\varepsilon,f}\|_\varepsilon^2
 +|c_{\varepsilon,f}|^2
 \leq\frac{C_q}{q}\|f\|_g^2
 \leq C_q\|f\|_g^2.
\]
This proves \eqref{geo-augmented-equation}.
\end{proof}

\subsection{Limits of solutions}

\begin{lemma}\label{geo-limit}
Let $\varepsilon_\nu\in(0,1]$ with
$\varepsilon_\nu\downarrow0$, and let
$u_\nu\in L^2_{n,q-1}(M,E;\omega_{\varepsilon_\nu},h)$ satisfy
$\sup_\nu\|u_\nu\|_{\varepsilon_\nu}<\infty$.
Suppose that
\[
 \bar\partial u_\nu=f_\nu
 \text{ in the sense of distribution},
 \qquad
 f_\nu\longrightarrow f
 \text{ in }L^2_{n,q}(M,E;\omega_g,h).
\]
Then there exists a subsequence  of  $\{u_\nu\}$ converging  to an $E$-valued $(n,q-1)$-form
$u$ weakly in $L_g^2(V)$ for each relatively compact open set
$V\Subset M$.
Moreover, $u\in\operatorname{Dom}(\bar\partial_g^{n,q-1})$ and
\[
 \bar\partial_g^{n,q-1}u=f,
 \qquad
 \|u\|_g\leq
 \liminf_{\nu\to\infty}\|u_\nu\|_{\varepsilon_\nu}.
\]
\end{lemma}

\begin{proof}

Choose a subsequence $(u_{\nu_k})$ such that
\begin{equation}\label{equa-lim-inf}
\lim_{k\to\infty}\|u_{\nu_k}\|_{\eps_{\nu_k}}
=
\liminf_{\nu\to\infty}\|u_\nu\|_{\eps_\nu}.
\end{equation}
Replace the original sequence by this subsequence
and keep the same notation.
Choose compact sets $K_j\subset M$ such that
\[
K_j\subset\operatorname{int}K_{j+1},
\qquad
\bigcup_{j=1}^\infty\operatorname{int}K_j=M.
\]

For each fixed $j$, choose $C_{K_j}\geq0$ such that
$\widehat\omega\leq C_{K_j}\omega_g$ on $K_j$.
Then
\[
 \omega_g\leq\omega_{\varepsilon_\nu}
 \leq(1+C_{K_j}\varepsilon_\nu)\omega_g
 \qquad\text{on }K_j.
\]
By the metric comparison argument used in the proof of Lemma
\ref{geo-comparison}, applied 
to the above inequalities, we obtain
\[
 \|u_\nu\|_{L_g^2(K_j)}
 \leq(1+C_{K_j}\varepsilon_\nu)^{(q-1)/2}
      \|u_\nu\|_{L_{\varepsilon_\nu}^2(K_j)}.
\]
Note that $(1+C_{K_j}\varepsilon_\nu)^{(q-1)/2}\to1$
as $\nu\to\infty$ for each fixed $j$.
In particular, $(u_\nu)$ is bounded in $L_g^2(K_j)$ for each  $j$.

By weak compactness and a diagonal subsequence, there is
a locally square-integrable form $u$ such that
\[
 u_\nu\rightharpoonup u
 \qquad\text{in }L_g^2(K_j)
\]
for each $j$.
Weak lower semicontinuity then gives
\[
 \begin{aligned}
 \|u\|_{L_g^2(K_j)}
 &\leq\liminf_{\nu\to\infty}
       \|u_\nu\|_{L_g^2(K_j)}\\
 &\leq\liminf_{\nu\to\infty}
       \left(
       (1+C_{K_j}\varepsilon_\nu)^{(q-1)/2}
       \|u_\nu\|_{\varepsilon_\nu}
       \right)\\
 &=\liminf_{\nu\to\infty}
       \|u_\nu\|_{\varepsilon_\nu}.
 \end{aligned}
\]
Letting $j\to\infty$ and using monotone convergence,
we obtain $u\in L^2_{n,q-1}(M,E;\omega_g,h)$.
Moreover, by \eqref{equa-lim-inf} and the above estimate,
we obtain the desired
\[
 \|u\|_g\leq
 \liminf_{\nu\to\infty}\|u_\nu\|_{\varepsilon_\nu}.
\]

Since the local weak $L^2$ convergence implies convergence in
distributions and 
 $f_\nu\to f$ in $L^2_{n,q}(M,E;\omega_g,h)$, we obtain
\[
 \bar\partial u
 =\lim_{\nu\to\infty}\bar\partial u_\nu
 =\lim_{\nu\to\infty}f_\nu
 =f
\]
in the sense of distributions.
Since $u\in L^2_{n,q-1}(M,E;\omega_g,h)$ and $f\in L^2_{n,q}(M,E;\omega_g,h)$, the definition of
the maximal realization gives
$u\in\operatorname{Dom}(\bar\partial_g^{n,q-1})$ and $\bar\partial_g^{n,q-1}u=f$.  This completes the proof of Lemma \ref{geo-limit}.
\end{proof}

\subsection{Proof of Theorem \ref{thm-main}}\label{ss-singular-metrics}

We are now ready to prove Theorem \ref{thm-main}.

\begin{proof}[Proof of Theorem \ref{thm-main}]

We first prove the theorem when $h$ is smooth, and then reduce the
singular metric case to the smooth case by comparing the corresponding
$L^2$ complexes.


\step{1}  The case where $h$ is smooth

Fix $1\leq q\leq n$.
We apply Theorem \ref{fa-cancellation} to the setting
\[
 \begin{aligned}
 X_0&=L^2_{n,q-1}(M,E;\omega_g,h),\qquad Z_0=\ker(\bar\partial_g^{n,q}),\\
 X_\varepsilon&=L^2_{n,q-1}(M,E;\omega_\varepsilon,h),\qquad
 Z_\varepsilon=\ker(\bar\partial_\varepsilon^{n,q}).
 \end{aligned}
\]
Define
\[
 T_0:\operatorname{Dom}(\bar\partial_g^{n,q-1})\subset L^2_{n,q-1}(M,E;\omega_g,h)
 \longrightarrow \ker(\bar\partial_g^{n,q}),
 \qquad
 T_0u=\bar\partial_g^{n,q-1}u.
\]
Then  $T_0$ is closed.
For $\varepsilon>0$, we use the operators $T_\varepsilon$
as defined in \eqref{geo-target}.

For $f\in \ker(\bar\partial_g^{n,q})$, Lemma \ref{geo-comparison} gives
$\|f\|_\varepsilon\leq\|f\|_g$, and hence $f\in \ker(\bar\partial_\varepsilon^{n,q})$.
Thus the linear map
\[
 I_\varepsilon:\ker(\bar\partial_g^{n,q})\longrightarrow \ker(\bar\partial_\varepsilon^{n,q}),
 \qquad I_\varepsilon f=f
\]
is well defined.
For each $f\in \ker(\bar\partial_g^{n,q})$, Lemma \ref{geo-augmented} gives
$u_{\varepsilon,f}\in\operatorname{Dom}T_\varepsilon$
and $c_{\varepsilon,f}\in\mathbb C^N$ such that
\[
 \begin{aligned}
&\qquad T_\varepsilon u_{\varepsilon,f}
 =I_\varepsilon f-L_\varepsilon^*c_{\varepsilon,f},\\
& \sup_{0<\varepsilon\leq1}
 \|u_{\varepsilon,f}\|_\varepsilon
 \leq\sqrt{C_q}\|f\|_g<\infty.
 \end{aligned}
\]
This verifies \eqref{fa-finite-error} with
$G_\varepsilon=-L_\varepsilon^*$.

Next, let $\varepsilon_\nu\downarrow0$,
$f_\nu,f\in \ker(\bar\partial_g^{n,q})$, and
$u_\nu\in\operatorname{Dom}T_{\varepsilon_\nu}$ satisfy
\[
 \|f_\nu-f\|_g\longrightarrow0,
 \qquad
 \sup_\nu\|u_\nu\|_{\varepsilon_\nu}<\infty,
 \qquad
 T_{\varepsilon_\nu}u_\nu=I_{\varepsilon_\nu}f_\nu.
\]
Then
Lemma \ref{geo-limit} gives
$u\in\operatorname{Dom}(\bar\partial_g^{n,q-1})$ such that
\[
 \bar\partial_g^{n,q-1}u=f.
\]
This verifies assumption \ref{fa-cancellation-ii}
of Theorem \ref{fa-cancellation}.
Applying Theorem \ref{fa-cancellation}, we obtain
\[
 \dim_{\mathbb C}(\ker(\bar\partial_g^{n,q})/\operatorname{Im}(\bar\partial_g^{n,q-1}))\leq N,
 \qquad
 \overline{\operatorname{Im}(\bar\partial_g^{n,q-1})}^{\,\ker(\bar\partial_g^{n,q})}=\operatorname{Im}(\bar\partial_g^{n,q-1}).
\]
Since $\ker(\bar\partial_g^{n,q})$ is closed in $L^2_{n,q}(M,E;\omega_g,h)$,
this proves Theorem \ref{thm-main} when $h$ is smooth.

\par\medskip

  \step{2}  The case where $h$ is singular

In this step, we construct a suitable smooth metric and use it to compare
the corresponding $L^2$ complexes with the $L^2$ complexes associated
with the singular metric $h$.

\par\medskip

\noindent \textbf{Step 2.1}  Construction of a smooth  metric

Set $S=\supp(\omega_g-\Theta_h)$. Note that on $M\setminus S$, the local
weight $\varphi$ satisfies
$\sqrt{-1}\,\partial\bar\partial\varphi=\omega_g$
in the sense of currents. Taking the Euclidean trace in local coordinates
then shows that $\Delta\varphi$ is smooth on $M\setminus S$.
By elliptic regularity, the weights have smooth representatives on
$M\setminus S$. Thus $h$ is smooth on $M\setminus S$.

Choose any smooth Hermitian metric $h_{\rm sm}$ on $E$.
On $M\setminus S$, write
\[
 h=e^{-\psi}h_{\rm sm}
\]
for a smooth real-valued function $\psi$. Choose compact sets
$K_0\subset\operatorname{int}K$ such that
$S\subset\operatorname{int}K_0$, and choose a smooth cutoff function
$\chi:M\to[0,1]$ such that $\chi=0$ near $K_0$ and
$\chi=1$ on $M\setminus K$. Extending $\chi\psi$ by zero across
$K_0$, define
\[
 h_0:=e^{-\chi\psi}h_{\rm sm}.
\]
Then $h_0$ is a smooth Hermitian metric on $E$ and
$h_0=h$ on $M\setminus K$.

The local weights of $h$ are quasi-plurisubharmonic, hence upper
semicontinuous and therefore bounded above on compact subsets of
their coordinate domains. Since $h_0$ is smooth, its local weights
are bounded below on compact subsets. Hence a finite cover of $K$
gives a constant $C\geq1$ such that
\begin{equation}\label{sing-metric-comparison}
 h_0\leq Ch\quad\text{on }M,
 \qquad
 h_0=h\quad\text{on }M\setminus K .
\end{equation}
Moreover, $\Theta_{h_0}=\omega_g$ on $M\setminus K$, and thus
$\supp(\omega_g-\Theta_{h_0})\subset K.$
Therefore results obtained in  Step~1 applies to $(E,h_0)$ with the original base metric $g$.

\par\medskip
\noindent\textbf{Step 2.2} The quotient of the local $L^2$ complexes

Recall the local $L^2$ complexes
$\mathcal L_h^{n,\bullet}$ from Subsection
\ref{ss-local-complexes}. Define
$\mathcal L_{h_0}^{n,\bullet}$ in the same way, with $h_0$
in place of $h$. By the metric comparison in Step 2.1, we have $\mathcal L_h^{n,\bullet}
 \subset
 \mathcal L_{h_0}^{n,\bullet}.$
Note that both complexes consist of fine sheaves. Moreover,
$\mathcal L_h^{n,\bullet}$ is a fine resolution of
$K_M\otimes E\otimes\mathcal I(h)$, while
$\mathcal L_{h_0}^{n,\bullet}$ is a fine resolution of
$K_M\otimes E$.

Put
\[
 \mathcal B^\bullet
 :=
 \mathcal L_{h_0}^{n,\bullet}/
 \mathcal L_h^{n,\bullet}.
\]
Since $h=h_0$ on $M\setminus K$, we have
\[
 \mathcal B^j|_{M\setminus K}=0.
\]

\par\medskip
\noindent\textbf{Step 2.3} The global  $L^2$ complexes and a short exact sequence of complexes

Let
\[
 D_h^\bullet:
 D_h^0\xrightarrow{\bar\partial}
 D_h^1\xrightarrow{\bar\partial}\cdots
 \xrightarrow{\bar\partial}D_h^n
\]
and
\[
 D_{h_0}^\bullet:
 D_{h_0}^0\xrightarrow{\bar\partial}
 D_{h_0}^1\xrightarrow{\bar\partial}\cdots
 \xrightarrow{\bar\partial}D_{h_0}^n
\]
be the complexes defined by
\[
 D_h^j
 =
 \left\{
 u\in L^2_{n,j}(M,E;\omega_g,h):
 \bar\partial u\in L^2_{n,j+1}(M,E;\omega_g,h)
 \right\},
\]
and
\[
 D_{h_0}^j
 =
 \left\{
 u\in L^2_{n,j}(M,E;\omega_g,h_0):
 \bar\partial u\in L^2_{n,j+1}(M,E;\omega_g,h_0)
 \right\},
\]
By definition, for each $1\leq q\leq n$, we have
\[
 H^q(D_h^\bullet)
 =
 H_{(2),\max}^{n,q}(M,E;g,h),
 \qquad
 H^q(D_{h_0}^\bullet)
 =
 H_{(2),\max}^{n,q}(M,E;g,h_0).
\]

We claim that the natural sequence
\begin{equation}\label{shor-exat-Bj}
     0\longrightarrow D_h^\bullet
 \longrightarrow D_{h_0}^\bullet
 \longrightarrow \Gamma(M,\mathcal B^\bullet)
 \longrightarrow0
\end{equation}
is a short exact sequence of complexes of complex vector spaces.

Indeed, for each $j$, the first map
$D_h^j\to D_{h_0}^j$ is the natural inclusion and is injective
by the metric comparison in Step 2.1. We next prove that
\[
 D_{h_0}^j\longrightarrow\Gamma(M,\mathcal B^j)
\]
is surjective. For each $j$, the definition of $\mathcal B^j$
gives a short exact sequence of sheaves
\[
 0\longrightarrow \mathcal L_h^{n,j}
 \longrightarrow \mathcal L_{h_0}^{n,j}
 \longrightarrow \mathcal B^j
 \longrightarrow0.
\]
Since $\mathcal L_h^{n,j}$ is fine, we have
$H^1(M,\mathcal L_h^{n,j})=0$. Therefore the induced map
\[
 \Gamma(M,\mathcal L_{h_0}^{n,j})
 \longrightarrow
 \Gamma(M,\mathcal B^j)
\]
is surjective, and thus any section
$s\in\Gamma(M,\mathcal B^j)$ lifts to a section
$\widetilde s\in\Gamma(M,\mathcal L_{h_0}^{n,j})$.
Since $\mathcal B^j|_{M\setminus K}=0$, choose
$\eta\in C_c^\infty(M)$ such that $\eta=1$ on a neighborhood
of $K$. Then $\eta\widetilde s$ has the same image as
$\widetilde s$ in $\Gamma(M,\mathcal B^j)$.
Moreover,
$\eta\widetilde s\in D_{h_0}^j$, which proves the surjectivity.

It remains to prove that for each $j$,
\[
 \ker\bigl(D_{h_0}^j\to\Gamma(M,\mathcal B^j)\bigr)=D_h^j.
\]
Clearly, $D_h^j$ maps to zero in $\Gamma(M,\mathcal B^j)$.
Conversely, if $u\in D_{h_0}^j$ maps to zero, then $u$ belongs
locally to $\mathcal L_h^{n,j}$. A finite relatively compact
cover of $K$ therefore shows that both $u$ and $\bar\partial u$
are square integrable over $K$ with respect to $h$.
On $M\setminus K$, the same conclusion follows from
$h=h_0$ and $u\in D_{h_0}^j$. Thus $u\in D_h^j$.

Finally, the natural inclusion and quotient maps commute with
the distributional $\bar\partial$. Therefore \eqref{shor-exat-Bj}
is a short exact sequence of complexes of complex vector spaces.

\par\medskip
\noindent\textbf{Step 2.4} Finiteness and closed range

By Step 2.2, taking the termwise quotient of the natural inclusion
$\mathcal L_h^{n,\bullet}\subset\mathcal L_{h_0}^{n,\bullet}$ of fine resolutions, we obtain an exact sequence
\[
 0\longrightarrow
 \frac{K_M\otimes E}
      {K_M\otimes E\otimes\mathcal I(h)}
 \longrightarrow
 \mathcal B^0
 \longrightarrow\cdots\longrightarrow
 \mathcal B^n
 \longrightarrow0.
\]
Denoted by 
\[\mathcal Q:=
 (K_M\otimes E)\otimes_{\mathcal O_M}
 \bigl(\mathcal O_M/\mathcal I(h)\bigr)\cong
 \frac{K_M\otimes E}
      {K_M\otimes E\otimes\mathcal I(h)}.
\]
Then $\mathcal B^\bullet$ is a fine resolution of $\mathcal Q$.
Consequently,  for each $j\geq0$, we have
\begin{equation}\label{iso-Hj-B-Hj}
     H^j\bigl(\Gamma(M,\mathcal B^\bullet)\bigr)
 \cong H^j(M,\mathcal Q).
\end{equation}

 Let $Z:=V(\mathcal I(h))$ be the zero variety, and
let $i:Z\hookrightarrow M$ be the inclusion. Since
$\mathcal I(h)=\mathcal O_M$ on $M\setminus K$, we have
$Z\subset K$, and thus $Z$ is compact. Moreover,
\[
 \mathcal Q
 \cong
 i_*\bigl((K_M\otimes E)|_Z\bigr),
\]
and therefore
\[
 H^j(M,\mathcal Q)
 \cong
 H^j\bigl(Z,(K_M\otimes E)|_Z\bigr)
\]
for each $j\geq0$. Since $Z$ is  compact (not necessarily reduced),
$ H^j\bigl(Z,(K_M\otimes E)|_Z\bigr)$ is finite dimensional by the Cartan--Serre
finiteness theorem (e.g., \cite[Chapter IX, Theorem 4.8]{Dem12}).

Applying the long exact cohomology sequence (e.g., \cite[Chapter IV, \S 1, (1.6)]{Dem12}) associated with
\eqref{shor-exat-Bj}, we obtain the exact segment
\[
 H^{q-1}\bigl(\Gamma(M,\mathcal B^\bullet)\bigr)
 \longrightarrow
 H_{(2),\max}^{n,q}(M,E;g,h)
 \longrightarrow
 H_{(2),\max}^{n,q}(M,E;g,h_0).
\]
By \eqref{iso-Hj-B-Hj}, we then obtain
 the exact segment
\[
 H^{q-1}(M,\mathcal Q)
 \longrightarrow
 H_{(2),\max}^{n,q}(M,E;g,h)
 \longrightarrow
 H_{(2),\max}^{n,q}(M,E;g,h_0).
\]
The LHS term is finite dimensional by the last
paragraph, while the RHS term is finite dimensional by results obtained in 
Step 1. This proves the desired
\[
 \dim_{\mathbb C}
 H_{(2),\max}^{n,q}(M,E;g,h)<\infty.
\]

Finally, by Fact \ref{fact-closed-den},
$\bar\partial_g^{n,q-1}$ is closed, and thus an argument similar to that at the end of
the proof of Theorem \ref{fa-cancellation} shows that
$\operatorname{Im}(\bar\partial_g^{n,q-1})$ is closed in
$L^2_{n,q}(M,E;\omega_g,h)$.  This completes the proof of Theorem  \ref{thm-main}.
\end{proof}

\section{Proof  of   Theorem \ref{thm-ohsawa-question-two-singular-intro}}\label{sec-ohsawa-question-two}

\subsection{An answer to Question \ref{ques-ohsawa-intro-2}}\label{ss-ohsawa-question-two-smooth}

In this subsection, we give an affirmative answer to Ohsawa's
\cite[Question 2]{Oh22} and will extend the answer to the singular
metric setting in the next subsection \ref{sec-ohsawa-question-two-singular}.

For establishing Theorem \ref{thm-ohsawa-question-two}, except for the construction of the complete approximating metrics,
where we use a plurisubharmonic exhaustion of the weakly $1$-complete
manifold, the proof is almost identical to that  of
Theorem \ref{thm-main} for smooth $h$.

\begin{theorem}[{Answer to Question \ref{ques-ohsawa-intro-2}}]\label{thm-ohsawa-question-two}
Let $M$ be a connected weakly $1$-complete complex manifold of
dimension $n\geq1$, and let $(E,h)$ be a smooth Hermitian
holomorphic line bundle on $M$ with curvature form $\Theta_h$.
Let $g$ be a smooth Hermitian metric on $M$ with fundamental form
$\omega_g$. Assume that
\[
 \supp(\omega_g-\Theta_h)\Subset M.
\]
Then for every $1\leq q\leq n$, the $L^2$ cohomology group
$H_{(2),\max}^{n,q}(M,E;g,h)$ in Definition
\ref{def-ltwo-cohomology} is finite dimensional.
In particular, $\operatorname{Im}(\bar\partial_g^{n,q-1})$
is closed in $L^2_{n,q}(M,E;\omega_g,h)$.
\end{theorem}

\begin{proof}
We follow the proof of Theorem~\ref{thm-main} for smooth $h$, using a
plurisubharmonic exhaustion of $M$ to construct the approximating metrics.
Fix $1\leq q\leq n$. 

\step{1} Complete metric approximation

Let $\rho$ be a smooth plurisubharmonic exhaustion function of $M$.
Since the  sublevel sets of  $\rho$  are relatively compact, $\rho$ is
bounded below. After adding a constant, we may assume that
$\rho\geq0$. Motivated by Ohsawa's (\cite[Proposition~1]{Oh81})  and Demailly's (\cite[proof of Theorem~1.3]{Dem82})
constructions,
 for $0\leq\eps\leq1$, we set
\begin{equation}\label{geo-weak-regularization}
 \omega_\eps:=\omega_g+
       \eps\sqrt{-1}\,\partial\bar\partial(\rho^2).
\end{equation}
The identity
\[
 \sqrt{-1}\,\partial\bar\partial(\rho^2)
 =2\rho\sqrt{-1}\,\partial\bar\partial\rho
   +2\sqrt{-1}\,\partial\rho\wedge\bar\partial\rho
\]
shows that the perturbation $\eps\sqrt{-1}\,\partial\bar\partial(\rho^2)$ is semipositive. For $\eps>0$,
the second term gives
$|d\rho|_{\omega_\eps}^2=2|\partial\rho|_{\omega_\eps}^2
\leq\eps^{-1}$.
Thus $\rho$ is Lipschitz for the
distance induced by $\omega_\eps$. Every Cauchy sequence
for this distance lies in a compact sublevel set of $\rho$ and
therefore has a convergent subsequence, so $\omega_\eps$ is complete for $\eps>0$.

Note that on every relatively compact open set $U\Subset M$, the
perturbing form is bounded above by $C_U\omega_g$ for some
$C_U\geq0$. It follows that
\begin{equation}\label{geo-weak-local-comparison}
 \omega_g\leq\omega_\eps\leq(1+C_U\eps)\omega_g
 \qquad\text{on }\overline U.
\end{equation}
The family depends smoothly on $\eps$, with
$\omega_0=\omega_g$.

In the rest of the proof, we use the notation of Section
\ref{sec-geometry} with \eqref{geo-weak-regularization} in place
of \eqref{geo-regularization}. Note that   the rest of the proof is almost identical to that  of
Theorem \ref{thm-main} for smooth $h$.

\step{2} Uniform estimate with a cutoff term

Define $A_{\eps,q}$ and $T_\eps$ as in
\eqref{Aepsi-aepsi} and \eqref{geo-target}, respectively. 
Note that the proof of Lemma \ref{geo-comparison}
uses only the pointwise inequality $\omega_\eps\geq\omega_g$.
Thus \eqref{Aepsilon-qI} and \eqref{geo-datum} hold for the
family \eqref{geo-weak-regularization} in our present setting.

Set $B=\Theta_h-\omega_g$. Since both the added form in
\eqref{geo-weak-regularization} and $\Theta_h$ are closed, we have
\[
 \partial\omega_\eps=-\partial B,\qquad
 \partial\bar\partial\omega_\eps=-\partial\bar\partial B.
\]
In  \eqref{geo-bkn},
the curvature term on $(n,q)$-forms is
$A_{\eps,q}+[L_B,\Lambda_\eps]$.  Moreover,  the error $[L_B,\Lambda_\eps]$
and the torsion term $\mathcal T_{\omega_\eps}$  in \eqref{geo-torsion} are supported in
$\supp B$. 

Note that the pointwise operator norms of $[L_B,\Lambda_\eps]$
and $\mathcal T_{\omega_\eps}$ are uniformly bounded on
$\supp B$ by the compactness of $\supp B\times[0,1]$
and the smooth dependence of the coefficients on $(x,\eps)$.
Choosing a fixed  cutoff $\chi\in C_c^\infty(M)$ with
$0\leq\chi\leq1$ and $\chi=1$ near $\supp B$, we obtain
\eqref{geo-basic-estimate} on compactly supported smooth
forms by discarding the nonnegative term in \eqref{geo-bkn}. 
For each $\eps>0$, by Lemma \ref{lem-graph-core} and
\eqref{Aepsilon-qI}, we pass the estimate
\eqref{geo-basic-estimate} to all
$v\in\Dom(\bar\partial_\varepsilon^{n,q})
\cap\Dom((\bar\partial_\varepsilon^{n,q-1})^*)$.
By \eqref{geo-adjoint}, we then obtain
\eqref{geo-closed-estimate} for $v\in\Dom T_\eps^*$.

\step{3} Uniform estimate with finitely many  linear functionals

As in the proof of Proposition \ref{geo-moments},
we replace the cutoff term by finitely many  linear
functionals.
Choose a fixed relatively compact open neighborhood $U$
of $\supp\chi$ as in the proof of Proposition
\ref{geo-moments}.
By the calculation in the proof of Lemma
\ref{geo-comparison} and \eqref{geo-weak-local-comparison},
we obtain \eqref{geo-local-uniform} on $\overline U$ in the present setting.

By the smooth dependence on $\eps$, we obtain
the uniform coefficient bounds needed in the  elliptic estimate
\eqref{geo-local-garding}. Hence by the same argument, we obtain  the estimate
\eqref{geo-sobolev} for the present family.
The Rellich compactness argument in Proposition \ref{geo-moments} now
applies in the fixed local $L^2$ space defined by $g$.
We then obtain the estimates \eqref{geo-moment-estimate} and
\eqref{fucn-esti-j1l} with $N$ and the functionals
$\ell_1,\ldots,\ell_N$ independent of $\eps$ (in the present setting).

\step{4} Application of the  criterion for closed range of finite codimension

We now apply Theorem \ref{fa-cancellation} to the setting
\[
 X_\eps=L^2_{n,q-1}(M,E;\omega_\eps,h),\qquad
 Z_\eps=\ker(\bar\partial_\eps^{n,q}),\qquad 0\leq\eps\leq1,
\]
and with $T_\eps$ defined as in \eqref{geo-target}
for each $0\leq\eps\leq1$.  Note that these operators $T_\eps$ are
closed with the indicated targets.
By \eqref{geo-datum} and the definition of the maximal
realization, the inclusion $I_\eps f=f$ from
$Z_0$ to $Z_\eps$ is well defined.
Let $L_\eps:Z_\eps\to\C^N$ be the bounded
map with components $\ell_1,\ldots,\ell_N$, as in the proof
of Lemma \ref{geo-augmented}.
As in the proof of Lemma \ref{geo-augmented},
\eqref{geo-moment-estimate} and \eqref{geo-datum}, together with
Lemma \ref{fa-augmented}, give \eqref{geo-augmented-equation}
for the present family.
This verifies condition  \eqref{fa-finite-error} of Theorem \ref{fa-cancellation} with $G_\eps=-L_\eps^*$.

To verify condition \ref{fa-cancellation-ii} of Theorem \ref{fa-cancellation}, suppose that
$\eps_\nu\downarrow0$, $f_\nu\to f$ in $Z_0$, and $u_\nu$
satisfies \eqref{fa-exact-limit-hypothesis}. The local comparison of the $L^2$ norms (used in the proof of
Lemma \ref{geo-limit}) follows in the present setting
from \eqref{geo-weak-local-comparison}.
In fact,
\[
 \|u_\nu\|_{L_g^2(U)}
 \leq(1+C_U\eps_\nu)^{(q-1)/2}\|u_\nu\|_{\eps_\nu}.
\]
The factor tends to $1$ on each fixed $U\Subset M$.
As in the proof of Lemma \ref{geo-limit}, a diagonal argument gives
$u\in X_0$ with $\bar\partial u=f$ in the sense of distributions.
Thus $u\in\Dom T_0$ and $T_0u=f$ by Definition
\ref{maximal-realization}.

It then follows from Theorem \ref{fa-cancellation}  that
$H_{(2),\max}^{n,q}(M,E;g,h)=Z_0/\operatorname{Im}T_0$ is finite
dimensional and that $\operatorname{Im}T_0$ is closed in $Z_0$.
Since $Z_0$ is closed in $L^2_{n,q}(M,E;\omega_g,h)$,
$\operatorname{Im}T_0$ is closed in this space as well.  This completes the proof of Theorem \ref{thm-ohsawa-question-two}.
\end{proof}

\subsection{Generalization to the singular metric case}\label{sec-ohsawa-question-two-singular}
We now extend Theorem \ref{thm-ohsawa-question-two}
to singular Hermitian metrics with quasi-plurisubharmonic local
weights.

\begin{theorem}[{= Theorem \ref{thm-ohsawa-question-two-singular-intro}}]
\label{thm-ohsawa-question-two-singular}
Let $M$ be a connected weakly $1$-complete (not necessarily K\"ahler) complex manifold of
dimension $n\geq1$. Let $E$ be a holomorphic line bundle
on $M$ with a singular Hermitian metric $h$ whose local weights
are quasi-plurisubharmonic. Let $g$ be a smooth Hermitian metric
on $M$ with the fundamental form $\omega_g$. Assume that 
\[
 \supp(\omega_g-\Theta_h)\Subset M,
\]
where the support is in the  sense of current.
Then for every $1\leq q\leq n$,
the $L^2$ cohomology group $H_{(2),\max}^{n,q}(M,E;g,h)$
in Definition \ref{def-ltwo-cohomology} is finite dimensional.
In particular, $\operatorname{Im}(\bar\partial_g^{n,q-1})$
is closed in $L^2_{n,q}(M,E;\omega_g,h)$.
\end{theorem}

\begin{proof}
The proof is the same as the singular metric part of the proof of
Theorem \ref{thm-main}. Indeed, by the same argument, there exist a
smooth Hermitian metric $h_0$ on $E$, a compact subset $K\Subset M$,
and a constant $C\geq1$ such that
\[
 h_0\leq Ch\quad\text{on }M,
 \qquad
 h_0=h\quad\text{on }M\setminus K,
\]
and
\[
 \supp(\omega_g-\Theta_{h_0})\subset K.
\]
Thus Theorem \ref{thm-ohsawa-question-two} applies to $h_0$.
The comparison of the local and global $L^2$ complexes in
Steps 2.2--2.4 of the proof of Theorem \ref{thm-main} then applies
without any change, and gives the desired finite dimensionality and
closed range property.
\end{proof}

\noindent\textbf{Acknowledgements.}
The authors would like to express their sincere gratitude to Professors Xiaoshan Li and Sheng Rao for  many valuable discussions regarding $L^2$ and Sobolev estimates, as well as other related topics.

\par\medskip
\par\medskip

\noindent\textbf{AI Declaration.}
The authors used ChatGPT as an auxiliary tool to deal with  the singular
metric case of Theorem \ref{thm-main}.
   The authors also used DeepSeek to revise sentences and improve the organization and readability. After using these tools, the authors reviewed and edited the content as needed and take full responsibility for the content of the paper.

\end{document}